\documentclass[12pt]{amsart}
\usepackage[colorlinks=true,citecolor=black,linkcolor=black,urlcolor=blue]{hyperref}
\usepackage{amsmath}
\usepackage[english,  activeacute]{babel}
\usepackage[utf8]{inputenc}
\usepackage{amssymb}
\usepackage{amsthm}
\usepackage{graphics,graphicx}
\usepackage{enumerate}
\usepackage{ytableau}
\usepackage{array}
\usepackage{bm}
\usepackage{cite}
\usepackage{a4wide}
\usepackage{color, url}
\usepackage{float}
\usepackage{comment}
\usepackage{xcolor}
\usepackage{color, url}
\usepackage{float}
\usepackage{pgf, tikz}

\theoremstyle{plain}
\newtheorem{theorem}{Theorem}[section]

\newtheorem{corollary}[theorem]{Corollary}

\theoremstyle{definition}

\newcommand{\area}{{\texttt{area}}}
\newcommand{\inter}{{\texttt{inter}}}
\newcommand{\sper}{{\texttt{sper}}}

\usepackage[usestackEOL]{stackengine}[2013-10-15]
\stackMath

\title{Enumeration on polyominoes determined by Catalan words avoiding $(\geq,\geq)$}

\date{\today}

\begin{document}
\author[M. Ahmia]{Moussa Ahmia}
\address{LMAM laboratory, University of Mohamed Seddik Benyahia, BP 98 Ouled Aissa 18000 Jijel, Algeria}
\email{moussa.ahmia@univ-jijel.dz;
ahmiamoussa@gmail.com}

\author[J-L. Baril]{Jean-Luc Baril}
\address{LIB, Université de Bourgogne Franche-Comté, B.P. 47 870, 21078, Dijon Cedex, France}
\email{barjl@u-bourgogne.fr}

\author[B. Rezig]{Boualam Rezig}
\address{Department of Informatic and Mathematics, High Normal School Constantine \\  LMAM laboratory, BP 98 Ouled Aissa 18000 Jijel, Algeria}
\email{rezig.boualam@ensc.dz; boualem.rezig@gmail.com}

\newcommand{\nadji}[1]{\mbox{}{\sf\color{green}[Kessouri: #1]}\marginpar{\color{green}\Large$*$}} 

\newcommand{\ahmia}[1]{\mbox{}{\sf\color{magenta}[Ahmia: #1]}\marginpar{\color{magenta}\Large$*$}} 

\newcommand{\jluc}[1]{\mbox{}{\sf\color{blue}[jluc: #1]}\marginpar{\color{blue}\Large$*$}}

\begin{abstract}
A Catalan word of length $n$ that avoids the pattern $(\geq, \geq)$ is a sequence 
$w=w_1\cdots w_n$ with $w_1=0$ and $0\leq w_i\leq w_{i-1}+1$  for all $i$, while ensuring that no subsequence satisfies $w_i \geq w_{i+1}\geq w_{i+2}$ for $i=2,\ldots,n$. These words are enumerated by the $n$-th Motzkin number. From such a word, we associate a $n$-column Motzkin polyomino (called a $(\geq,\geq)$-polyomino), where the $i$-th column contains $w_i+1$ bottom-aligned cells. In this paper, we derive generating functions for $(\geq,\geq)$-polyominoes based on their length, area, semiperimeter, last symbol value, and number of interior points. We provide asymptotic analyses and closed-form expressions for the total area, total semiperimeter, sum of the last symbol values, and total number of interior points across all $(\geq,\geq)$-polyominoes of a given length. Finally, we express all these results as linear combinations of  trinomial coefficients.
\end{abstract}

\maketitle
\section{Introduction}
\begin{figure}[htb]
    \centering
   \scalebox{0.7}{ \includegraphics[width=0.9\linewidth]{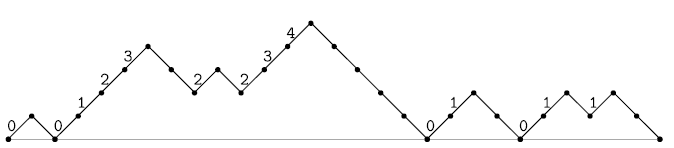}}
        \centering
   \scalebox{0.4}{\begin{tikzpicture}
% Définir la couleur de remplissage
%\definecolor{lightpink}{rgb}{1.0, 0.71, 0.76}
\definecolor{lightpink}{rgb}{1.0, 0.85, 0.88}
% Dessiner les cellules du polyomino avec leurs délimitations
% Colonne 1 (1 cellule)
\fill[lightpink] (0,0) rectangle (1,1);
\draw[black] (0,0) rectangle (1,1);
% Colonne 2 (2 cellules)
\fill[lightpink] (1,0) rectangle (2,1);
%\fill[lightpink] (1,1) rectangle (2,2);
\draw[black] (1,0) rectangle (2,1);
%\draw[black] (1,1) rectangle (2,2);
% Colonne 3 (3 cellules)
\fill[lightpink] (2,0) rectangle (3,1);
\fill[lightpink] (2,1) rectangle (3,2);
\draw[black] (2,0) rectangle (3,1);
\draw[black] (2,1) rectangle (3,2);
% Colonne 4 (1 cellule)
\fill[lightpink] (3,0) rectangle (4,1);
\fill[lightpink] (3,1) rectangle (4,2);
\fill[lightpink] (3,2) rectangle (4,3);
\draw[black] (3,0) rectangle (4,1);
\draw[black] (3,1) rectangle (4,2);
\draw[black] (3,2) rectangle (4,3);

% Colonne 5 (1 cellule)
\fill[lightpink] (4,0) rectangle (5,1);
\fill[lightpink] (4,1) rectangle (5,2);
\fill[lightpink] (4,2) rectangle (5,3);
\fill[lightpink] (4,3) rectangle (5,4);
\draw[black] (4,0) rectangle (5,1);
\draw[black] (4,1) rectangle (5,2);
\draw[black] (4,2) rectangle (5,3);
\draw[black] (4,3) rectangle (5,4);
% Colonne 5 (1 cellule)
%\fill[lightpink] (4,0) rectangle (5,1);
%\draw[black] (4,0) rectangle (5,1);

% Colonne 6 (1 cellule)
\fill[lightpink] (5,0) rectangle (6,1);
\fill[lightpink] (5,1) rectangle (6,2);
\fill[lightpink] (5,2) rectangle (6,3);
\draw[black] (5,0) rectangle (6,1);
\draw[black] (5,1) rectangle (6,2);
\draw[black] (5,2) rectangle (6,3);

% Colonne 7 (1 cellule)
\fill[lightpink] (6,0) rectangle (7,1);
\fill[lightpink] (6,1) rectangle (7,2);
\fill[lightpink] (6,2) rectangle (7,3);
\draw[black] (6,0) rectangle (7,1);
\draw[black] (6,1) rectangle (7,2);
\draw[black] (6,2) rectangle (7,3);

% Colonne 8 (1 cellule)
\fill[lightpink] (7,0) rectangle (8,1);
\fill[lightpink] (7,1) rectangle (8,2);
\fill[lightpink] (7,2) rectangle (8,3);
\fill[lightpink] (7,3) rectangle (8,4);
\draw[black] (7,0) rectangle (8,1);
\draw[black] (7,1) rectangle (8,2);
\draw[black] (7,2) rectangle (8,3);
\draw[black] (7,3) rectangle (8,4);

% Colonne 9 (1 cellule)
\fill[lightpink] (8,0) rectangle (9,1);
\fill[lightpink] (8,1) rectangle (9,2);
\fill[lightpink] (8,2) rectangle (9,3);
\fill[lightpink] (8,3) rectangle (9,4);
\fill[lightpink] (8,4) rectangle (9,5);
\draw[black] (8,0) rectangle (9,1);
\draw[black] (8,1) rectangle (9,2);
\draw[black] (8,2) rectangle (9,3);
\draw[black] (8,3) rectangle (9,4);
\draw[black] (8,4) rectangle (9,5);
% Colonne 10 (1 cellule)
\fill[lightpink] (9,0) rectangle (10,1);
\draw[black] (9,0) rectangle (10,1);

% Colonne 11 (1 cellule)
\fill[lightpink] (10,0) rectangle (11,1);
\fill[lightpink] (10,1) rectangle (11,2);
\draw[black] (10,0) rectangle (11,1);
\draw[black] (10,1) rectangle (11,2);

% Colonne 12 (1 cellule)
\fill[lightpink] (11,0) rectangle (12,1);
\draw[black] (11,0) rectangle (12,1);

% Colonne 13 (1 cellule)
\fill[lightpink] (12,0) rectangle (13,1);
\fill[lightpink] (12,1) rectangle (13,2);
\draw[black] (12,0) rectangle (13,1);
\draw[black] (12,1) rectangle (13,2);

% Colonne 14 (1 cellule)
\fill[lightpink] (13,0) rectangle (14,1);
\fill[lightpink] (13,1) rectangle (14,2);
\draw[black] (13,0) rectangle (14,1);
\draw[black] (13,1) rectangle (14,2);
%\node at (6,-0.35){\Large \texttt{011201123011}};
\end{tikzpicture}} 
    \caption{Dyck path associated with the Catalan word $\texttt{00123223401011}$, and its corresponding polyomino.}
    \label{fig1}
\end{figure}
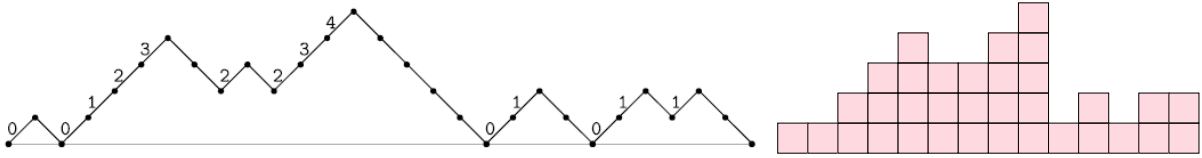

A {\it Catalan word} of length $n\geq 0$ is a sequence $w=w_1w_2\cdots w_n$ over the set of non-negative integers that satisfies $w_1=0$ and $0 \leq w_i\leq w_{i-1}+1$ for $i=2,\ldots,n$. Let $\mathcal{C}_n$ be the set of all Catalan words of length $n$, and let $\epsilon$ be the empty word (i.e. the Catalan word of length $0$). The number of such words is given by the $n$-th {\it Catalan number} $C_n=\frac{1}{n+1}\binom{2n}{n}$, see \cite[Exercise 80]{stan}. There is a very close link between Catalan words and  some paths lying into the first quadrant of the plane. Indeed, from {\it Dyck paths} of semilength $n$, i.e. lattice paths of $\mathbb{N}^2$ starting at $(0,0)$, ending at $(2n,0)$ and consisting of up-steps $U=(1, 1)$ and down-steps $D=(1, -1)$, we define a Catalan word in $\mathcal{C}_n$ by recording  the $y$-coordinates of the up step starting points (read from left to right).  For instance, the left part of Figure~\ref{fig1} illustrates the Dyck path associated with the Catalan word $\texttt{00123223401011}\in \mathcal{C}_{14}$. See for instance \cite{Banfla,Ges,Krat} for references about lattice paths.

A Catalan word $w = w_1 \cdots w_n$ can also be viewed as a polyomino (also called a bargraph) where the $i$-th column contains $w_i+1$ cells for $1 \leq i \leq n$, and all columns are bottom-justified (i.e., the bottom-most cells of all columns are aligned in the same row). The polyomino associated with a Catalan word of length $n$ is called a \emph{Catalan polyomino} of length $n$. In Figure~\ref{fig1}, we show the Catalan polyomino of length 14 associated with the Catalan word \texttt{00123223401011}.

Let $w$ be a Catalan word and $P(w)$ its associated polyomino. We denote by $\area(w)$ the number of cells (or the \emph{area}) of $P(w)$, which is also the sum of all $1+w_i$ for $1 \leq i \leq n$. The \emph{semiperimeter} of $P(w)$, denoted $\sper(w)$, is half of the \emph{perimeter} of $P(w)$, which is the number of cell borders that do not touch another cell of $P(w)$. An \emph{interior point} of $P(w)$ is a point that belongs to exactly four cells of $P(w)$. We denote by $\inter(w)$ the number of interior points of $P(w)$. For instance, if $w$ is the polyomino in Figure~\ref{fig1}, then $\area(w) = 34$, $\sper(w) = 22$, and $\inter(w) = 13$. We refer to \cite{Book1} for a historical review on polyominoes, and to \cite{ManSha2} for definitions of many statistical and enumerative methods related to polyominoes. Additionally, analogous results are known for bargraphs \cite{Bou, BLE3}, compositions \cite{BleBreKnop}, set partitions \cite{ManA2, ManA3}, Catalan words \cite{CallManRam, Toc, ManRam}, inversion sequences \cite{ArcBleKnop}, and words \cite{ManA, BleBreKnop3}.

In \cite{BKR},  Catalan polyominoes of length $n$ that do not have two adjacent columns of the same height are studied.  In our paper, these polyominoes are called {\it $(\neq)$-polyominoes}. Mansour and Ram\'irez (see \cite{ManRamMot}) proved that they are counted by the Motzkin number $m_{n-1}$, where $m_n$ is defined by the  combinatorial sum (see \cite{RiordanN, Don})
$$m_n=\frac{1}{n+1}\sum_{i\geq 0}\binom{n+1}{i}\binom{n+1-i}{i+1}, \quad n\geq 0.$$
The number $m_n$  corresponds to the $n$-th term of the generating function 
$$M(x):=\sum_{n\geq 0}m_nx^n=\frac{1-x-\sqrt{1-2x-3x^2}}{2x^2}.$$ In \cite{BKR}, the authors give enumerative results for the number of these polyominoes according to the area, the  semiperimeter and the number of interior points, and they show the  surprising fact that the number of cells of a given height in all polyominoes of length $n$ corresponds to the trinomial coefficients. Recently, a comparable result for Catalan polyominoes was given by  Blecher and Knopfmacher \cite{BlKn} who prove a relationship between the number of cells at different heights and the first terms of the expanded polynomial $(1+x)^{2n}$.
\medskip

\begin{figure}[htb]
\centering
\scalebox{0.4}{\begin{tikzpicture}
% Définir la couleur de remplissage
%\definecolor{lightpink}{rgb}{1.0, 0.71, 0.76}
\definecolor{lightpink}{rgb}{1.0, 0.85, 0.88}
% Dessiner les cellules du polyomino avec leurs délimitations
% Colonne 1 (1 cellule)
\fill[lightpink] (0,0) rectangle (1,1);
\draw[black] (0,0) rectangle (1,1);
% Colonne 2 (2 cellules)
\fill[lightpink] (1,0) rectangle (2,1);
%\fill[lightpink] (1,1) rectangle (2,2);
\draw[black] (1,0) rectangle (2,1);
%\draw[black] (1,1) rectangle (2,2);
% Colonne 3 (3 cellules)
\fill[lightpink] (2,0) rectangle (3,1);
\fill[lightpink] (2,1) rectangle (3,2);
%\fill[lightpink] (2,2) rectangle (3,3);
\draw[black] (2,0) rectangle (3,1);
\draw[black] (2,1) rectangle (3,2);
%\draw[black] (2,2) rectangle (3,3);
% Colonne 4 (1 cellule)
\fill[lightpink] (3,0) rectangle (4,1);
\draw[black] (3,0) rectangle (4,1);
\node at (2,-0.45){\Large \texttt{0~0~1~0}};
\end{tikzpicture}}%*****************
\quad 
\scalebox{0.4}{\begin{tikzpicture}
% Définir la couleur de remplissage
%\definecolor{lightpink}{rgb}{1.0, 0.71, 0.76}
\definecolor{lightpink}{rgb}{1.0, 0.85, 0.88}
% Dessiner les cellules du polyomino avec leurs délimitations
% Colonne 1 (1 cellule)
\fill[lightpink] (0,0) rectangle (1,1);
\draw[black] (0,0) rectangle (1,1);
% Colonne 2 (2 cellules)
\fill[lightpink] (1,0) rectangle (2,1);
%\fill[lightpink] (1,1) rectangle (2,2);
\draw[black] (1,0) rectangle (2,1);
%\draw[black] (1,1) rectangle (2,2);
% Colonne 3 (3 cellules)
\fill[lightpink] (2,0) rectangle (3,1);
\fill[lightpink] (2,1) rectangle (3,2);
%\fill[lightpink] (2,2) rectangle (3,3);
\draw[black] (2,0) rectangle (3,1);
\draw[black] (2,1) rectangle (3,2);
%\draw[black] (2,2) rectangle (3,3);
% Colonne 4 (1 cellule)
\fill[lightpink] (3,0) rectangle (4,1);
\fill[lightpink] (3,1) rectangle (4,2);
\draw[black] (3,0) rectangle (4,1);
\draw[black] (3,1) rectangle (4,2);
\node at (2,-0.45){\Large \texttt{0~0~1~1}};\end{tikzpicture}}
\quad 
\scalebox{0.4}{\begin{tikzpicture}
% Définir la couleur de remplissage
%\definecolor{lightpink}{rgb}{1.0, 0.71, 0.76}
\definecolor{lightpink}{rgb}{1.0, 0.85, 0.88}
% Dessiner les cellules du polyomino avec leurs délimitations
% Colonne 1 (1 cellule)
\fill[lightpink] (0,0) rectangle (1,1);
\draw[black] (0,0) rectangle (1,1);
% Colonne 2 (2 cellules)
\fill[lightpink] (1,0) rectangle (2,1);
%\fill[lightpink] (1,1) rectangle (2,2);
\draw[black] (1,0) rectangle (2,1);
%\draw[black] (1,1) rectangle (2,2);
% Colonne 3 (3 cellules)
\fill[lightpink] (2,0) rectangle (3,1);
\fill[lightpink] (2,1) rectangle (3,2);
%\fill[lightpink] (2,2) rectangle (3,3);
\draw[black] (2,0) rectangle (3,1);
\draw[black] (2,1) rectangle (3,2);
%\draw[black] (2,2) rectangle (3,3);
% Colonne 4 (1 cellule)
\fill[lightpink] (3,0) rectangle (4,1);
\fill[lightpink] (3,1) rectangle (4,2);
\fill[lightpink] (3,2) rectangle (4,3);
\draw[black] (3,0) rectangle (4,1);
\draw[black] (3,1) rectangle (4,2);
\draw[black] (3,2) rectangle (4,3);
\node at (2,-0.45){\Large \texttt{0~0~1~2}};
\end{tikzpicture}}
\quad 
\scalebox{0.4}{\begin{tikzpicture}
% Définir la couleur de remplissage
%\definecolor{lightpink}{rgb}{1.0, 0.71, 0.76}
\definecolor{lightpink}{rgb}{1.0, 0.85, 0.88}
% Dessiner les cellules du polyomino avec leurs délimitations
% Colonne 1 (1 cellule)
\fill[lightpink] (0,0) rectangle (1,1);
\draw[black] (0,0) rectangle (1,1);
% Colonne 2 (2 cellules)
\fill[lightpink] (1,0) rectangle (2,1);
\fill[lightpink] (1,1) rectangle (2,2);
\draw[black] (1,0) rectangle (2,1);
\draw[black] (1,1) rectangle (2,2);
% Colonne 3 (3 cellules)
\fill[lightpink] (2,0) rectangle (3,1);
%\fill[lightpink] (2,1) rectangle (3,2);
%\fill[lightpink] (2,2) rectangle (3,3);
\draw[black] (2,0) rectangle (3,1);
%\draw[black] (2,1) rectangle (3,2);
%\draw[black] (2,2) rectangle (3,3);
% Colonne 4 (1 cellule)
\fill[lightpink] (3,0) rectangle (4,1);
\fill[lightpink] (3,1) rectangle (4,2);
\draw[black] (3,0) rectangle (4,1);
\draw[black] (3,1) rectangle (4,2);
\node at (2,-0.45){\Large \texttt{0~1~0~1}};
\end{tikzpicture}}
\quad 
\scalebox{0.4}{\begin{tikzpicture}
% Définir la couleur de remplissage
%\definecolor{lightpink}{rgb}{1.0, 0.71, 0.76}
\definecolor{lightpink}{rgb}{1.0, 0.85, 0.88}
% Dessiner les cellules du polyomino avec leurs délimitations
% Colonne 1 (1 cellule)
\fill[lightpink] (0,0) rectangle (1,1);
\draw[black] (0,0) rectangle (1,1);
% Colonne 2 (2 cellules)
\fill[lightpink] (1,0) rectangle (2,1);
\fill[lightpink] (1,1) rectangle (2,2);
\draw[black] (1,0) rectangle (2,1);
\draw[black] (1,1) rectangle (2,2);
% Colonne 3 (3 cellules)
\fill[lightpink] (2,0) rectangle (3,1);
\fill[lightpink] (2,1) rectangle (3,2);
%\fill[lightpink] (2,2) rectangle (3,3);
\draw[black] (2,0) rectangle (3,1);
\draw[black] (2,1) rectangle (3,2);
%\draw[black] (2,2) rectangle (3,3);
% Colonne 4 (1 cellule)
\fill[lightpink] (3,0) rectangle (4,1);
\fill[lightpink] (3,1) rectangle (4,2);
\fill[lightpink] (3,2) rectangle (4,3);
\draw[black] (3,0) rectangle (4,1);
\draw[black] (3,1) rectangle (4,2);
\draw[black] (3,2) rectangle (4,3);
\node at (2,-0.45){\Large \texttt{0~1~1~2}};
\end{tikzpicture}}
\quad 
\scalebox{0.4}{\begin{tikzpicture}
% Définir la couleur de remplissage
%\definecolor{lightpink}{rgb}{1.0, 0.71, 0.76}
\definecolor{lightpink}{rgb}{1.0, 0.85, 0.88}
% Dessiner les cellules du polyomino avec leurs délimitations
% Colonne 1 (1 cellule)
\fill[lightpink] (0,0) rectangle (1,1);
\draw[black] (0,0) rectangle (1,1);
% Colonne 2 (2 cellules)
\fill[lightpink] (1,0) rectangle (2,1);
\fill[lightpink] (1,1) rectangle (2,2);
\draw[black] (1,0) rectangle (2,1);
\draw[black] (1,1) rectangle (2,2);
% Colonne 3 (3 cellules)
\fill[lightpink] (2,0) rectangle (3,1);
\fill[lightpink] (2,1) rectangle (3,2);
\fill[lightpink] (2,2) rectangle (3,3);
\draw[black] (2,0) rectangle (3,1);
\draw[black] (2,1) rectangle (3,2);
\draw[black] (2,2) rectangle (3,3);
% Colonne 4 (1 cellule)
\fill[lightpink] (3,0) rectangle (4,1);
%\fill[lightpink] (3,1) rectangle (4,2);
%\fill[lightpink] (3,2) rectangle (4,3);
\draw[black] (3,0) rectangle (4,1);
%\draw[black] (3,1) rectangle (4,2);
%\draw[black] (3,2) rectangle (4,3);
\node at (2,-0.45){\Large \texttt{0~1~2~0}};
\end{tikzpicture}}
\quad \scalebox{0.4}{\begin{tikzpicture}
% Définir la couleur de remplissage
%\definecolor{lightpink}{rgb}{1.0, 0.71, 0.76}
\definecolor{lightpink}{rgb}{1.0, 0.85, 0.88}
% Dessiner les cellules du polyomino avec leurs délimitations
% Colonne 1 (1 cellule)
\fill[lightpink] (0,0) rectangle (1,1);
\draw[black] (0,0) rectangle (1,1);
% Colonne 2 (2 cellules)
\fill[lightpink] (1,0) rectangle (2,1);
\fill[lightpink] (1,1) rectangle (2,2);
\draw[black] (1,0) rectangle (2,1);
\draw[black] (1,1) rectangle (2,2);
% Colonne 3 (3 cellules)
\fill[lightpink] (2,0) rectangle (3,1);
\fill[lightpink] (2,1) rectangle (3,2);
\fill[lightpink] (2,2) rectangle (3,3);
\draw[black] (2,0) rectangle (3,1);
\draw[black] (2,1) rectangle (3,2);
\draw[black] (2,2) rectangle (3,3);
% Colonne 4 (1 cellule)
\fill[lightpink] (3,0) rectangle (4,1);
\fill[lightpink] (3,1) rectangle (4,2);
%\fill[lightpink] (3,2) rectangle (4,3);
\draw[black] (3,0) rectangle (4,1);
\draw[black] (3,1) rectangle (4,2);
%\draw[black] (3,2) rectangle (4,3);
\node at (2,-0.45){\Large \texttt{0~1~2~1}};
\end{tikzpicture}}
\quad \scalebox{0.4}{\begin{tikzpicture}
% Définir la couleur de remplissage
%\definecolor{lightpink}{rgb}{1.0, 0.71, 0.76}
\definecolor{lightpink}{rgb}{1.0, 0.85, 0.88}
% Dessiner les cellules du polyomino avec leurs délimitations
% Colonne 1 (1 cellule)
\fill[lightpink] (0,0) rectangle (1,1);
\draw[black] (0,0) rectangle (1,1);
% Colonne 2 (2 cellules)
\fill[lightpink] (1,0) rectangle (2,1);
\fill[lightpink] (1,1) rectangle (2,2);
\draw[black] (1,0) rectangle (2,1);
\draw[black] (1,1) rectangle (2,2);
% Colonne 3 (3 cellules)
\fill[lightpink] (2,0) rectangle (3,1);
\fill[lightpink] (2,1) rectangle (3,2);
\fill[lightpink] (2,2) rectangle (3,3);
\draw[black] (2,0) rectangle (3,1);
\draw[black] (2,1) rectangle (3,2);
\draw[black] (2,2) rectangle (3,3);
% Colonne 4 (1 cellule)
\fill[lightpink] (3,0) rectangle (4,1);
\fill[lightpink] (3,1) rectangle (4,2);
\fill[lightpink] (3,2) rectangle (4,3);
\draw[black] (3,0) rectangle (4,1);
\draw[black] (3,1) rectangle (4,2);
\draw[black] (3,2) rectangle (4,3);
\node at (2,-0.45){\Large \texttt{0~1~2~2}};
\end{tikzpicture}}
\quad \scalebox{0.4}{\begin{tikzpicture}
% Définir la couleur de remplissage
%\definecolor{lightpink}{rgb}{1.0, 0.71, 0.76}
\definecolor{lightpink}{rgb}{1.0, 0.85, 0.88}
% Dessiner les cellules du polyomino avec leurs délimitations
% Colonne 1 (1 cellule)
\fill[lightpink] (0,0) rectangle (1,1);
\draw[black] (0,0) rectangle (1,1);
% Colonne 2 (2 cellules)
\fill[lightpink] (1,0) rectangle (2,1);
\fill[lightpink] (1,1) rectangle (2,2);
\draw[black] (1,0) rectangle (2,1);
\draw[black] (1,1) rectangle (2,2);
% Colonne 3 (3 cellules)
\fill[lightpink] (2,0) rectangle (3,1);
\fill[lightpink] (2,1) rectangle (3,2);
\fill[lightpink] (2,2) rectangle (3,3);
\draw[black] (2,0) rectangle (3,1);
\draw[black] (2,1) rectangle (3,2);
\draw[black] (2,2) rectangle (3,3);
% Colonne 4 (1 cellule)
\fill[lightpink] (3,0) rectangle (4,1);
\fill[lightpink] (3,1) rectangle (4,2);
\fill[lightpink] (3,2) rectangle (4,3);
\fill[lightpink] (3,3) rectangle (4,4);
\draw[black] (3,0) rectangle (4,1);
\draw[black] (3,1) rectangle (4,2);
\draw[black] (3,2) rectangle (4,3);
\draw[black] (3,3) rectangle (4,4);
\node at (2,-0.45){\Large \texttt{0~1~2~3}};
\end{tikzpicture}}
\caption{All $(\geq, \geq)$-polyominoes of length $4$. }\label{fig2}
\end{figure}
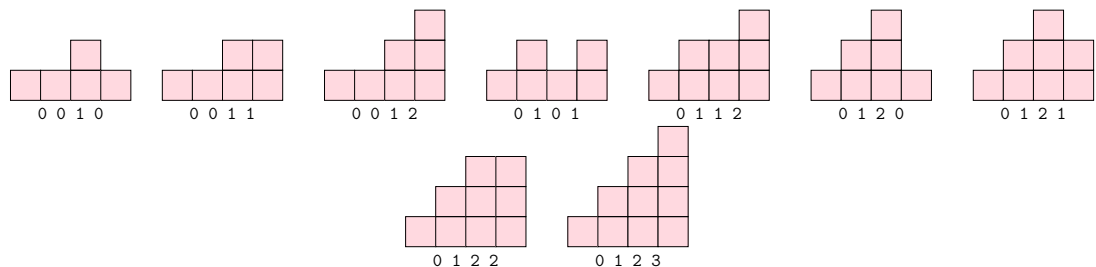

In our work, we consider the set $\mathcal{C}_{(\geq, \geq)}(n)$ of all Catalan words of length $n$ avoiding the pattern $(\geq, \geq)$, i.e. all Catalan words $w$ of length $n$ such that there is no $i$ satisfying $w_i\geq w_{i+1}\geq w_{i+2}$.   The avoidance of the pattern $(\geq, \geq)$ on Catalan words is equivalent to the avoidance of the consecutive patterns $\underline{000},\underline{100}, \underline{110}$ and $\underline{210}$. For example,
\[\mathcal{C}_{(\geq, \geq)}(4)=\{\texttt{0010},\texttt{0011},\texttt{0012},\texttt{0101},\texttt{0112},\texttt{0120},\texttt{0121},\texttt{0122},\texttt{0123}\}.
\]
Polyominoes associated with Catalan words in $\mathcal{C}_{(\geq, \geq)}(n)$ will be called {\it $(\geq, \geq)$-polyominoes}.
 Figure \ref{fig2} illustrates all $(\geq, \geq)$-polyominoes of length $4$. In \cite{relation}, it is proved that the cardinality of $\mathcal{C}_{(\geq, \geq)}(n)$ is given by the $n$-th Motzkin number $m_n$ defined above.

\medskip

 \noindent{\bf Outline of the paper.} The motivation for this work comes from the two recent studies \cite{BKR,BlKn} on Catalan polyominoes. We provide enumerative results on polyominoes associated to Catalan words avoiding $(\geq,\geq)$ according to several parameters (area, semiperimeter, and number of interior points). The results are obtained algebraically by considering functional equations for multivariate  generating functions.
In Section \ref{sec:2}, we focus on the statistics of the area, the semiperimeter and the last symbol. We provide the functional equation for multivariate generating function, where the coefficient of $x^np^k q^\ell v^r$ is the number of $(\geq,\geq)$-polyominoes $P(w)$ of length $n$, satisfying $\sper(w)=k$ and $\area(w)=\ell$ and $w_n=r$.
From this result, we deduce the generating function for the total sum of the last symbol over all $(\geq,\geq)$-avoiding Catalan words of length $n$, and we give asymptotic approximation for this number. 
In Section \ref{sec:3}, we focus on the semiperimeter statistic. We give the generating function for nonempty $(\geq,\geq)$-polyominoes with respect to the length and the semiperimeter.
We also exhibit a bijection between $(\geq,\geq)$-polyominoes of length $n$  and ($\neq$)-polyominoes  of length $n+1$ that send a semiperimeter $k$ into a semiperimeter $k+2$. We provide a closed-form expression for the total semiperimeter over all $(\geq,\geq)$-polyominoes of length $n$.
In Section \ref{sec:4}, we make a similar study as in Section \ref{sec:3} by considering the area instead of the semiperimeter. Finally, Section \ref{sec:5} is dedicated to the statistic of the number of interior points of the $(\geq,\geq)$-polyominoes, and we exhibit a closed-form expression for the total number of interior points over all $(\geq,\geq)$-polyominoes of length $n$. As in \cite{BKR,BlKn}, our study make a link between the number of cells over all $(\geq,\geq)$-polyominoes of length $n$ and the $n$-th trinomial coefficient, i.e. the coefficient of $x^n$ in the expanded polynomial $(1+x+x^2)^n$.

\section{Area, semiperimeter and last symbol statistics}\label{sec:2}
In this section, we investigate the area, the semiperimeter and the last symbol statistics on the class of polyominoes associated with Catalan words avoiding  the pattern $(\geq, \geq)$. For $0\leq i \leq n-1$, let $\mathcal{C}_{(\geq, \geq)}(n,i)$ denote the set of such Catalan words of length $n$ whose last symbol is $i$. We  introduce the generating functions
\[
C_i^{\geq}(x;p,q)=\sum_{n\geq 1}x^n \sum_{w\in \mathcal{C}_{(\geq, \geq)}(n,i)}p^{\sper(w)}q^{\area(w)},
\]
and 
\[
C^{\geq}(x;p,q;v)=\sum_{i \geq 0}C_i^{\geq}(x;p,q)v^{i}.
\]
\begin{theorem}
The generating function $C^{\geq}(x;p,q;v)$ satisfies the functional equation
\begin{align}
C^{\geq}(x;p,q;v)=p^2qx+p^3q^2x^2+\frac{p^3q^3x^2}{1-qv}C^{\geq}(x;p,q;q)&+p^2 q^2xvC^{\geq}(x;p,q;qv)\\ \notag
&-\frac{p^3q^5 x^2v^2}{1-qv}C^{\geq}(x;p,q;q^2v).
\end{align}\label{feq}
\end{theorem}
\begin{proof}
Let $w$ be a Catalan word in $\mathcal{C}_{(\geq, \geq)}(n)$. We analyze two distinct cases:  (i) $w\in \mathcal{C}_{(\geq, \geq)}(n,0)$, and (ii) $w\in \mathcal{C}_{(\geq, \geq)}(n,i)$ for $i\geq 1$. Remember that if $w\in \mathcal{C}_{(\geq, \geq)}(n,i)$ then the corresponding polyomino $P(w)$ contains $i+1$ cells in its rightmost column.

Case (i). The last symbol of $w$ is $0$. Hence, we have either $w=0$ or $w=00$ or $w=w'k0$ where $w'\in \mathcal{C}_{(\geq, \geq)}(n-2,k-1)$, with $n\geq 3$ and $k>0$. See Figure \ref{f14} for a graphical representation of these three cases. 

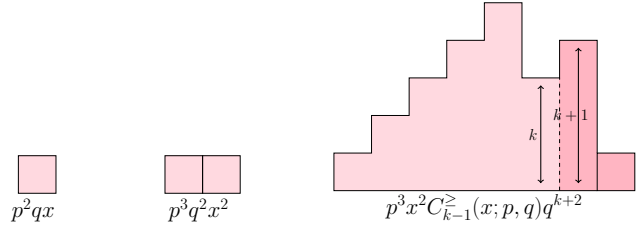
\begin{figure}[htb]
\centering

\scalebox{0.5}{\begin{tikzpicture}
% Définir la couleur de remplissage
%\definecolor{lightpink}{rgb}{1.0, 0.71, 0.76}
\definecolor{lightpink}{rgb}{1.0, 0.85, 0.88}
% Dessiner les cellules du polyomino avec leurs délimitations
% Colonne 1 (1 cellule)
\fill[lightpink] (0,0) rectangle (1,1);
\draw[black] (0,0) rectangle (1,1);
\node at (0.4,-0.45){\Large $p^2qx$};
\end{tikzpicture}}
\quad\quad\quad
\scalebox{0.5}{\begin{tikzpicture}
% Définir la couleur de remplissage
%\definecolor{lightpink}{rgb}{1.0, 0.71, 0.76}
\definecolor{lightpink}{rgb}{1.0, 0.85, 0.88}
% Dessiner les cellules du polyomino avec leurs délimitations
% Colonne 1 (1 cellule)
\fill[lightpink] (0,0) rectangle (1,1);
\fill[lightpink] (1,0) rectangle (2,1);
\draw[black] (0,0) rectangle (1,1);
\draw[black] (1,0) rectangle (2,1);
\node at (0.9,-0.45){\Large $p^3q^2x^2$};
\end{tikzpicture}}\quad\quad\quad
\scalebox{0.5}{\begin{tikzpicture}
% Définir la couleur de remplissage
\definecolor{morepink}{rgb}{1.0, 0.71, 0.76}
\definecolor{lightpink}{rgb}{1.0, 0.85, 0.88}
% Dessiner les cellules du polyomino avec leurs délimitations
% Colonne 1 (1 cellule)
\fill[lightpink] (0,0) rectangle (1,1);
\fill[lightpink] (1,0) rectangle (2,1);
\fill[lightpink] (1,1) rectangle (2,2);
\fill[lightpink] (2,0) rectangle (3,1);
\fill[lightpink] (2,1) rectangle (3,2);
\fill[lightpink] (2,2) rectangle (3,3);
\fill[lightpink] (3,0) rectangle (4,1);
\fill[lightpink] (3,1) rectangle (4,2);
\fill[lightpink] (3,2) rectangle (4,3);
\fill[lightpink] (3,3) rectangle (4,4);
\fill[lightpink] (4,0) rectangle (5,1);
\fill[lightpink] (4,1) rectangle (5,2);
\fill[lightpink] (4,2) rectangle (5,3);
\fill[lightpink] (4,3) rectangle (5,4);
\fill[lightpink] (4,4) rectangle (5,5);
\fill[lightpink] (5,0) rectangle (6,1);
\fill[lightpink] (5,1) rectangle (6,2);
\fill[lightpink] (5,2) rectangle (6,3);
\fill[morepink] (6,0) rectangle (7,1);
\fill[morepink] (6,1) rectangle (7,2);
\fill[morepink] (6,2) rectangle (7,3);
\fill[morepink] (6,3) rectangle (7,4);
\fill[morepink] (7,0) rectangle (8,1);
\draw[black] (7,0) rectangle (8,1);
\draw[black] (7,0)--(0,0)--(0,1)--(1,1)--(1,2)--(2,2)--(2,3)--(3,3)--(3,4)--(4,4)--(4,5)--(5,5)--(5,3)--(6,3)--(6,4)--(7,4)--(7,1);
\draw[<->, thick] (6.5,0.2) -- (6.5,3.8); 
\node at (6.3,2){$k+1$};
\draw[<->, thick] (5.5,0.2) -- (5.5,2.8); 
\node at (5.3,1.5){ $k$};
\draw[dashed] (6,0) -- (6,3); 
% Colonne 4 (1 cellule)
\node at (4,-0.45){\Large $p^3x^2C_{k-1}^{\geq}(x;p,q)q^{k+2}$};
\end{tikzpicture}}
 \caption{Case (i): The last column is of height $1$.}\label{f14}
 \end{figure}

The generating function for this case is given by

\[
C_0^{\geq}(x;p,q)=p^2 qx+p^3 q^2x^2+p^3x^2\sum_{k \geq 1}C_{k-1}^{\geq}(x;p,q)q^{k+2}.
 \]
This relation can also be written
 \begin{equation}\label{gf1}
C_0^{\geq}(x;p,q)=p^2 qx+p^3 q^2x^2+p^3q^3x^2C^{\geq}(x;p,q;q).
 \end{equation}

Case (ii). The last symbol of $w$ is at least $1$. From the definition of a word $w\in \mathcal{C}_{(\geq, \geq)}(n)$, we have the decompositions either $w=w'i$ with $w'\in \mathcal{C}_{(\geq, \geq)}(n-1,i-1) $, $n\geq 2$, or  $w=w'ki$, with $w'\in \mathcal{C}_{(\geq, \geq)}(n-2,k-1)$ and  $k \geq i$. See Figure \ref{f24} for a graphical representation of these two cases.

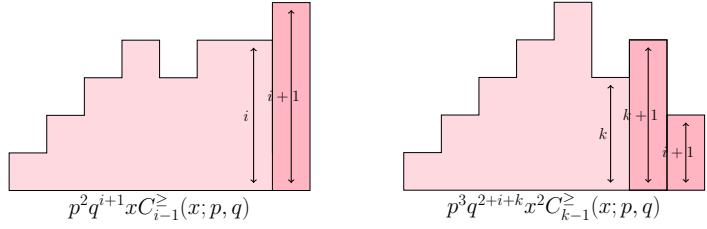
\begin{figure}[htb]
\centering

\scalebox{0.5}{\begin{tikzpicture}
% Définir la couleur de remplissage
\definecolor{morepink}{rgb}{1.0, 0.71, 0.76}
\definecolor{lightpink}{rgb}{1.0, 0.85, 0.88}
% Dessiner les cellules du polyomino avec leurs délimitations
% Colonne 1 (1 cellule)
\fill[lightpink] (0,0) rectangle (1,1);
\fill[lightpink] (1,0) rectangle (2,1);
\fill[lightpink] (1,1) rectangle (2,2);
\fill[lightpink] (2,0) rectangle (3,1);
\fill[lightpink] (2,1) rectangle (3,2);
\fill[lightpink] (2,2) rectangle (3,3);
\fill[lightpink] (3,0) rectangle (4,1);
\fill[lightpink] (3,1) rectangle (4,2);
\fill[lightpink] (3,2) rectangle (4,3);
\fill[lightpink] (3,3) rectangle (4,4);
\fill[lightpink] (4,0) rectangle (5,1);
\fill[lightpink] (4,1) rectangle (5,2);
\fill[lightpink] (4,2) rectangle (5,3);
%\fill[lightpink] (4,3) rectangle (5,4);
%\fill[lightpink] (4,4) rectangle (5,3);
\fill[lightpink] (5,0) rectangle (6,1);
\fill[lightpink] (5,1) rectangle (6,2);
\fill[lightpink] (5,2) rectangle (6,3);
\fill[lightpink] (5,3) rectangle (6,4);
\fill[lightpink] (6,0) rectangle (7,1);
\fill[lightpink] (6,1) rectangle (7,2);
\fill[lightpink] (6,2) rectangle (7,3);
\fill[lightpink] (6,3) rectangle (7,4);
\fill[morepink] (7,0) rectangle (8,5);
\draw[black] (7,0) rectangle (8,5);
\draw[black] (7,0)--(0,0)--(0,1)--(1,1)--(1,2)--(2,2)--(2,3)--(3,3)--(3,4)--(4,4)--(4,3)--(5,3)--(5,4)--(6,4)--(7,4)--(7,1);
\draw[<->, thick] (6.5,0.2) -- (6.5,3.8); 
\node at (6.3,2){ $i$};
\draw[<->, thick] (7.5,0.2) -- (7.5,4.8); 
\node at (7.3,2.5){ $i+1$};
% Colonne 4 (1 cellule)
\node at (4,-0.45){\Large $p^2q^{i+1}xC_{i-1}^{\geq}(x;p,q)$};
\end{tikzpicture}}\quad\quad\quad
\scalebox{0.5}{\begin{tikzpicture}
% Définir la couleur de remplissage
\definecolor{morepink}{rgb}{1.0, 0.71, 0.76}
\definecolor{lightpink}{rgb}{1.0, 0.85, 0.88}
% Dessiner les cellules du polyomino avec leurs délimitations
% Colonne 1 (1 cellule)
\fill[lightpink] (0,0) rectangle (1,1);
\fill[lightpink] (1,0) rectangle (2,1);
\fill[lightpink] (1,1) rectangle (2,2);
\fill[lightpink] (2,0) rectangle (3,1);
\fill[lightpink] (2,1) rectangle (3,2);
\fill[lightpink] (2,2) rectangle (3,3);
\fill[lightpink] (3,0) rectangle (4,1);
\fill[lightpink] (3,1) rectangle (4,2);
\fill[lightpink] (3,2) rectangle (4,3);
\fill[lightpink] (3,3) rectangle (4,4);
\fill[lightpink] (4,0) rectangle (5,1);
\fill[lightpink] (4,1) rectangle (5,2);
\fill[lightpink] (4,2) rectangle (5,3);
\fill[lightpink] (4,3) rectangle (5,4);
\fill[lightpink] (4,4) rectangle (5,5);
\fill[lightpink] (5,0) rectangle (6,1);
\fill[lightpink] (5,1) rectangle (6,2);
\fill[lightpink] (5,2) rectangle (6,3);
\fill[lightpink] (6,0) rectangle (7,1);
\fill[lightpink] (6,1) rectangle (7,2);
\fill[lightpink] (6,2) rectangle (7,3);
\fill[lightpink] (6,3) rectangle (7,4);
\fill[morepink] (7,0) rectangle (8,2);
\fill[morepink] (6,0) rectangle (7,4);
\draw[black] (6,0) rectangle (7,4);
\draw[black] (7,0) rectangle (8,2);
\draw[black] (7,0)--(0,0)--(0,1)--(1,1)--(1,2)--(2,2)--(2,3)--(3,3)--(3,4)--(4,4)--(4,5)--(5,5)--(5,3)--(6,3)--(6,4)--(7,4)--(7,1);
\draw[<->, thick] (6.5,0.2) -- (6.5,3.8); 
\node at (6.3,2){ $k+1$};
\draw[<->, thick] (5.5,0.2) -- (5.5,2.8); 
\node at (5.3,1.5){ $k$};
\draw[<->, thick] (7.5,0.2) -- (7.5,1.8); 
\node at (7.3,1){ $i+1$};

% Colonne 4 (1 cellule)
\node at (4,-0.45){\Large $p^3q^{2+i+k}x^2C_{k-1}^{\geq}(x;p,q)$};
\end{tikzpicture}}
 \caption{Case (ii): The last column is of height at least $2$.}\label{f24}
 \end{figure}

In terms of generating functions, we deduce the functional equations for any $i\geq 1$
 \begin{align}\label{gf2}
C_i^{\geq}(x;p,q)=p^2q^{i+1}xC_{i-1}^{\geq}(x;p,q)&+p^3q^{2i+2}x^2C_{i-1}^{\geq}(x;p,q) \\ \notag
& +p^3q^{i+2}x^2\sum_{k> i}q^kC_{k-1}^{\geq}(x;p,q).
 \end{align}

By multiplying Eq.~\eqref{gf2} by $v^i$ and summing over $i\geq1$, we obtain
\begin{align*}
C^{\geq}(x;p,q;v)-C_0^{\geq}(x;p,q)&= p^2q x\sum_{i\geq 1}v^iq^i C_{i-1}^{\geq}(x;p,q)+p^3 q^2x^2\sum_{i\geq 1}v^iq^{2i}C_{i-1}^{\geq}(x;p,q)\\
&\quad \quad \quad \quad\quad+p^3q^2 x^2\sum_{i\geq 1}v^iq^i \sum_{k> i}q^k C_{k-1}^{\geq}(x;p,q)\\
&= p^2q x\sum_{i\geq 1}v^iq^i C_{i-1}^{\geq}(x;p,q)+p^3 q^2x^2\sum_{i\geq 1}v^iq^{2i}C_{i-1}^{\geq}(x;p,q)\\
&\quad \quad \quad \quad\quad+p^3q^2 x^2\sum_{k\geq 2}q^k C_{k-1}^{\geq}(x;p,q)\sum_{i=1}^{k-1}v^iq^i\\
&= p^2 q^2xvC^{\geq}(x;p,q;qv)+p^3 q^4 x^2vC^{\geq}(x;p,q;q^2v)\\
&\quad \quad \quad \quad\quad+p^3q^2 x^2\sum_{k\geq 2}q^k C_{k-1}^{\geq}(x;p,q)\frac{qv-q^kv^k}{1-qv}\\
&=p^2 q^2xvC^{\geq}(x;p,q;qv)+p^3 q^4 x^2vC^{\geq}(x;p,q;q^2v)\\
&\quad \quad \quad+\frac{p^3q^4 x^2v}{1-qv}\left(C^{\geq}(x;p,q;q)-C_0^{\geq}(x;p,q)\right)\\
&\quad \quad \quad\quad\quad-\frac{p^3q^4 x^2v}{1-qv}\left(C^{\geq}(x;p,q;q^2v)-C_0^{\geq}(x;p,q)\right).
\end{align*}

Simplifying the last expression using Eq.~\eqref{gf1}, we obtain the desired result.
\end{proof}

Note that when $p=q=1$, we have
\[
\left(1-xv+\frac{x^2v^2}{1-v}\right)C^{\geq}(x;1,1;v)=x+x^2+\frac{x^2}{1-v}C^{\geq}(x;1,1;1).
\]
The generating function $C^{\geq}(x;1,1;1)$ corresponds to the generating function of the nonempty Catalan words in $ \mathcal{C}_{(\geq, \geq)}$, so $C^{\geq}(x;1,1;1)=M(x)-1$. Thus, we deduce the following.
\begin{corollary}
The generating function for the number of nonempty Catalan words with respect to the length and the value of
the last symbol is
\begin{equation}\label{gf3}
 C^{\geq}(x;1,1;v) =\frac{x(1-v)-x^2v+x^2M(x)}{1-v-xv(1-v)+x^2v^2}.
\end{equation}
\end{corollary}
The first terms of the series expansion are
\[
x+(v+1)x^2+(v^2+2v+1)x^3+(v^3+\mathbf{3v^2}+3v+2)x^4+(v^4+4v^3+6v^2+6v+4)x^5+O(x^6).
\]

There are three polyominoes of length $4$ ending with a column of height two (see Figure~\ref{fig1}). By differentiating $C^{\geq}(x;1,1;v) $
at $v = 1$, we deduce:

\begin{corollary}\label{corhn}
The generating function for the total sum $h(n)$ of the last symbol in all Catalan words of length $n$ avoiding the pattern $(\geq, \geq)$  is
\[\sum_{n\geq 1}h(n) x^n=\frac{\left(x +1\right) \left(\left(2x -1\right) \sqrt{-3 x^{2}-2 x +1}+2x^{3}-3 x+1\right)}{2x^{4}}.\] 
An asymptotic approximation for the coefficient $h(n)$ is given by 
$$\frac{2\sqrt{3}}{\sqrt{\pi}}\cdot\frac{3^{n+1}}{n^{\frac{3}{2}}}.$$
\end{corollary}

The asymptotic approximation is easily obtained with a singularity analysis using classical methods (see \cite{Fla,Orl}). The first terms of $h(n)$, $1\leq n\leq 10$, are $$0,\quad 1,\quad 4,\quad 12,\quad 34,\quad 94,\quad 258,\quad 707,\quad  1940, \quad 5337,$$ and the sequence $h(n)$ does not appear in \cite{OEIS}. Finally, the expected value of the last symbol is $h(n)/m_n$ and an asymptotic is $4$.

In the following corollary, we give an explicit relation for the sequence $h(n)$ using the central trinomial coefficient $T_n$, that is the coefficient of $x^n$ in the expanded polynomial $(1+x+x^2)^n$:
\[
T_n=\sum_{k=0}^{n}\binom{n}{k}\binom{n-k}{k}.
\]
Moreover, the generating function $T(x)$ of the central trinomial coefficients is given by
\begin{equation}\label{gft}
T(x):=\sum_{n\geq 0}T_n x^{n}=\frac{1}{\sqrt{1-2x-3x^2}}.
\end{equation}
\begin{corollary}The total sum of  the last symbol over all $(\geq,\geq)$-polyominoes of length $n$ is given by
$$h(n)=\frac{1}{2}\cdot\left(-6 T_n -7 T_{n+1}+3T_{n+2}+3T_{n+3}-T_{n+4}\right).$$
\end{corollary}
\begin{proof} We decompose the generating function of Corollary~\ref{corhn} as \[\frac{2 x^{4}+2 x^{3}-3 x^{2}-2 x +1}{2 x^{4}}+\frac{-6 x^{4}-7 x^{3}+3 x^{2}+3 x -1}{2 x^{4}}T(x).\]

Using relation \eqref{gft}, and comparing the $n$-th coefficient, we obtain the desired result.
\end{proof}

\section{The semiperimeter statistic}\label{sec:3}

In this section we study the semiperimeter statistic on $(\geq,\geq)$-polyominoes. By Eq.~\eqref{feq} with $q=1$, we obtain
\begin{align}
&\left(1-p^2 xv+\frac{p^3 x^2v^2}{1-v}\right)C^{\geq}(x;p,1;v)=p^2x+p^3x^2+\frac{p^3x^2}{1-v}C^{\geq}(x;p,1;1).\label{fq}
\end{align}
In order to compute $S(x, p) := C^{\geq}(x;p,1;1)$, we use the kernel method (see \cite{Kernel,pro}). The method consists in cancelling
the coefficient of $C^{\geq}(x;p,1;v)$ by taking the small root $v_0$ (the ones going to 0 for $x\thicksim 0$) of $$1-p^2 xv_0+\frac{p^3 x^2v_0^2}{1-v_0},$$ namely
\[
v_0=\frac{1+p^2x-\sqrt{1-2p^2x+(p^4-4p^3)x^2}}{2p^2x(1+px)}.
\]

When the factor on the left-hand side equals zero (and thus satisfies the equation), the right-hand side also equals zero
and so the solution $v_0$ can be used on the right-hand side to obtain $C^{\geq}(x;p,1;1)$ (see Theorem \ref{thm1}).
\begin{theorem}\label{thm1}The generating function for the number of nonempty $(\geq,\geq)$-polyominoes according to the length and the
semiperimeter is given by
\[
S(x, p)=\frac{1-p^2x-2p^3x^2-\sqrt{1-2p^2x+(p^4-4p^3)x^2}}{2p^3x^2}.
\]
\end{theorem}
The first terms of the series expansion of $S(x, p)$ are
\[
p^2x+(p^3+p^4)x^2+(3p^5+p^6)x^3+(2p^6+{\bf 6p^7}+p^8)x^4+(10p^8+10p^9+p^{10})x^5+O(x^6).
\]

\begin{figure}[htb]
\centering
\scalebox{0.4}{\begin{tikzpicture}
% Définir la couleur de remplissage
%\definecolor{lightpink}{rgb}{1.0, 0.71, 0.76}
\definecolor{lightpink}{rgb}{1.0, 0.85, 0.88}
% Dessiner les cellules du polyomino avec leurs délimitations
% Colonne 1 (1 cellule)
\fill[lightpink] (0,0) rectangle (1,1);
\draw[black] (0,0) rectangle (1,1);
% Colonne 2 (2 cellules)
\fill[lightpink] (1,0) rectangle (2,1);
%\fill[lightpink] (1,1) rectangle (2,2);
\draw[black] (1,0) rectangle (2,1);
%\draw[black] (1,1) rectangle (2,2);
% Colonne 3 (3 cellules)
\fill[lightpink] (2,0) rectangle (3,1);
\fill[lightpink] (2,1) rectangle (3,2);
%\fill[lightpink] (2,2) rectangle (3,3);
\draw[black] (2,0) rectangle (3,1);
\draw[black] (2,1) rectangle (3,2);
%\draw[black] (2,2) rectangle (3,3);
% Colonne 4 (1 cellule)
\fill[lightpink] (3,0) rectangle (4,1);
\fill[lightpink] (3,1) rectangle (4,2);
\fill[lightpink] (3,2) rectangle (4,3);
\draw[black] (3,0) rectangle (4,1);
\draw[black] (3,1) rectangle (4,2);
\draw[black] (3,2) rectangle (4,3);
\node at (2,-0.45){\Large 0~0~1~2};
\end{tikzpicture}}
\quad 
\scalebox{0.4}{\begin{tikzpicture}
% Définir la couleur de remplissage
%\definecolor{lightpink}{rgb}{1.0, 0.71, 0.76}
\definecolor{lightpink}{rgb}{1.0, 0.85, 0.88}
% Dessiner les cellules du polyomino avec leurs délimitations
% Colonne 1 (1 cellule)
\fill[lightpink] (0,0) rectangle (1,1);
\draw[black] (0,0) rectangle (1,1);
% Colonne 2 (2 cellules)
\fill[lightpink] (1,0) rectangle (2,1);
\fill[lightpink] (1,1) rectangle (2,2);
\draw[black] (1,0) rectangle (2,1);
\draw[black] (1,1) rectangle (2,2);
% Colonne 3 (3 cellules)
\fill[lightpink] (2,0) rectangle (3,1);
%\fill[lightpink] (2,1) rectangle (3,2);
%\fill[lightpink] (2,2) rectangle (3,3);
\draw[black] (2,0) rectangle (3,1);
%\draw[black] (2,1) rectangle (3,2);
%\draw[black] (2,2) rectangle (3,3);
% Colonne 4 (1 cellule)
\fill[lightpink] (3,0) rectangle (4,1);
\fill[lightpink] (3,1) rectangle (4,2);
\draw[black] (3,0) rectangle (4,1);
\draw[black] (3,1) rectangle (4,2);
\node at (2,-0.45){\Large 0~1~0~1};
\end{tikzpicture}}
\quad 
\scalebox{0.4}{\begin{tikzpicture}
% Définir la couleur de remplissage
%\definecolor{lightpink}{rgb}{1.0, 0.71, 0.76}
\definecolor{lightpink}{rgb}{1.0, 0.85, 0.88}
% Dessiner les cellules du polyomino avec leurs délimitations
% Colonne 1 (1 cellule)
\fill[lightpink] (0,0) rectangle (1,1);
\draw[black] (0,0) rectangle (1,1);
% Colonne 2 (2 cellules)
\fill[lightpink] (1,0) rectangle (2,1);
\fill[lightpink] (1,1) rectangle (2,2);
\draw[black] (1,0) rectangle (2,1);
\draw[black] (1,1) rectangle (2,2);
% Colonne 3 (3 cellules)
\fill[lightpink] (2,0) rectangle (3,1);
\fill[lightpink] (2,1) rectangle (3,2);
%\fill[lightpink] (2,2) rectangle (3,3);
\draw[black] (2,0) rectangle (3,1);
\draw[black] (2,1) rectangle (3,2);
%\draw[black] (2,2) rectangle (3,3);
% Colonne 4 (1 cellule)
\fill[lightpink] (3,0) rectangle (4,1);
\fill[lightpink] (3,1) rectangle (4,2);
\fill[lightpink] (3,2) rectangle (4,3);
\draw[black] (3,0) rectangle (4,1);
\draw[black] (3,1) rectangle (4,2);
\draw[black] (3,2) rectangle (4,3);
\node at (2,-0.45){\Large 0~1~1~2};
\end{tikzpicture}}
\quad 
\scalebox{0.4}{\begin{tikzpicture}
% Définir la couleur de remplissage
%\definecolor{lightpink}{rgb}{1.0, 0.71, 0.76}
\definecolor{lightpink}{rgb}{1.0, 0.85, 0.88}
% Dessiner les cellules du polyomino avec leurs délimitations
% Colonne 1 (1 cellule)
\fill[lightpink] (0,0) rectangle (1,1);
\draw[black] (0,0) rectangle (1,1);
% Colonne 2 (2 cellules)
\fill[lightpink] (1,0) rectangle (2,1);
\fill[lightpink] (1,1) rectangle (2,2);
\draw[black] (1,0) rectangle (2,1);
\draw[black] (1,1) rectangle (2,2);
% Colonne 3 (3 cellules)
\fill[lightpink] (2,0) rectangle (3,1);
\fill[lightpink] (2,1) rectangle (3,2);
\fill[lightpink] (2,2) rectangle (3,3);
\draw[black] (2,0) rectangle (3,1);
\draw[black] (2,1) rectangle (3,2);
\draw[black] (2,2) rectangle (3,3);
% Colonne 4 (1 cellule)
\fill[lightpink] (3,0) rectangle (4,1);
%\fill[lightpink] (3,1) rectangle (4,2);
%\fill[lightpink] (3,2) rectangle (4,3);
\draw[black] (3,0) rectangle (4,1);
%\draw[black] (3,1) rectangle (4,2);
%\draw[black] (3,2) rectangle (4,3);
\node at (2,-0.45){\Large 0~1~2~0};
\end{tikzpicture}}
\quad \scalebox{0.4}{\begin{tikzpicture}
% Définir la couleur de remplissage
%\definecolor{lightpink}{rgb}{1.0, 0.71, 0.76}
\definecolor{lightpink}{rgb}{1.0, 0.85, 0.88}
% Dessiner les cellules du polyomino avec leurs délimitations
% Colonne 1 (1 cellule)
\fill[lightpink] (0,0) rectangle (1,1);
\draw[black] (0,0) rectangle (1,1);
% Colonne 2 (2 cellules)
\fill[lightpink] (1,0) rectangle (2,1);
\fill[lightpink] (1,1) rectangle (2,2);
\draw[black] (1,0) rectangle (2,1);
\draw[black] (1,1) rectangle (2,2);
% Colonne 3 (3 cellules)
\fill[lightpink] (2,0) rectangle (3,1);
\fill[lightpink] (2,1) rectangle (3,2);
\fill[lightpink] (2,2) rectangle (3,3);
\draw[black] (2,0) rectangle (3,1);
\draw[black] (2,1) rectangle (3,2);
\draw[black] (2,2) rectangle (3,3);
% Colonne 4 (1 cellule)
\fill[lightpink] (3,0) rectangle (4,1);
\fill[lightpink] (3,1) rectangle (4,2);
%\fill[lightpink] (3,2) rectangle (4,3);
\draw[black] (3,0) rectangle (4,1);
\draw[black] (3,1) rectangle (4,2);
%\draw[black] (3,2) rectangle (4,3);
\node at (2,-0.45){\Large 0~1~2~1};
\end{tikzpicture}}
\quad \scalebox{0.4}{\begin{tikzpicture}
% Définir la couleur de remplissage
%\definecolor{lightpink}{rgb}{1.0, 0.71, 0.76}
\definecolor{lightpink}{rgb}{1.0, 0.85, 0.88}
% Dessiner les cellules du polyomino avec leurs délimitations
% Colonne 1 (1 cellule)
\fill[lightpink] (0,0) rectangle (1,1);
\draw[black] (0,0) rectangle (1,1);
% Colonne 2 (2 cellules)
\fill[lightpink] (1,0) rectangle (2,1);
\fill[lightpink] (1,1) rectangle (2,2);
\draw[black] (1,0) rectangle (2,1);
\draw[black] (1,1) rectangle (2,2);
% Colonne 3 (3 cellules)
\fill[lightpink] (2,0) rectangle (3,1);
\fill[lightpink] (2,1) rectangle (3,2);
\fill[lightpink] (2,2) rectangle (3,3);
\draw[black] (2,0) rectangle (3,1);
\draw[black] (2,1) rectangle (3,2);
\draw[black] (2,2) rectangle (3,3);
% Colonne 4 (1 cellule)
\fill[lightpink] (3,0) rectangle (4,1);
\fill[lightpink] (3,1) rectangle (4,2);
\fill[lightpink] (3,2) rectangle (4,3);
\draw[black] (3,0) rectangle (4,1);
\draw[black] (3,1) rectangle (4,2);
\draw[black] (3,2) rectangle (4,3);
\node at (2,-0.45){\Large 0~1~2~2};
\end{tikzpicture}}
\caption{The $6$ $(\geq,\geq)$-polyominoes of length $4$ and semiperimeter $7$. }\label{fig4}
\end{figure}
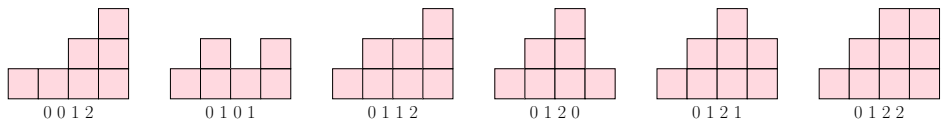

Figure~\ref{fig4} yields the  $6$ $(\geq,\geq)$-polyominoes of length $4$ and semiperimeter $7$.

 Notice that the number of $(\geq,\geq)$-polyominoes of length $n$ with a semiperimeter $k$ equals the number of $(\neq)$-Catalan polyominoes of length $n+1$ with a semiperimeter $k+2$ (see \cite{BKR}).

Let us give a bijection $\chi$ from the set $\mathcal{C}_{(\geq,\geq)}(n)$ to  the set $\mathcal{C}_{(\neq)}(n+1)$ of length $n+1$ Catalan words $w$ such that $w_i\neq w_{i+1}$ for $1\leq i \leq n-1$.

$$\chi(w)=\left\{\begin{array}{ll}
0 & \mbox{if } w=\epsilon, \\
 \texttt{0} (1+\chi(u))& \mbox{if } w=\texttt{0}(1+u),\\
  \texttt{0} (1+\chi(u))\texttt{0}& \mbox{if } w=\texttt{00}(1+u),\\
\texttt{0}(1+\chi(u')) \chi(v) & \mbox{if } w=\texttt{0}(1+u)v, \mbox{ where } u=u'a,v\neq \epsilon \mbox{ and } a\geq 1,
\end{array}\right.
$$
where $u'a$ and $v$ are $(\geq,\geq)$-Catalan words such that either $u'$ is empty or $u'$ ends with $a-1$, and $(1+u)$ corresponds to the word obtained from $u$ by increasing by one all letters of $u$. Using an induction on the length, we can easily check that the semiperimeter $k$ of a word $w$ is transported into the semiperimeter $k+2$ in the word $\chi(w)$. In order to convince the reader, let us assume the last case $w=0(1+u)v$ where $u=u'a$ and $v\neq \epsilon$, and let us prove by induction that $\sper(\chi(w))=\sper(w)+2$.
We have
\begin{align*}
\sper(\chi(w))=&\sper(\texttt{0}(1+\chi(u'))\chi(v))=\sper(\texttt{0}(1+\chi(u')))+\sper(\chi(v))-1\\
&=\sper(\chi(u'))-2+\sper(\chi(v))-1.
\end{align*}
Using the  recurrence hypothesis on $\chi(u')$ and $\chi(v)$, we obtain
\begin{align*}
\sper(\chi(w))&=\sper(u')+2-2+\sper(v)+2-1\\
&=\sper(u)+2+\sper(v)+1\\
&=\sper(\texttt{0}(1+u)v)+2=\sper(w)+2.
\end{align*}

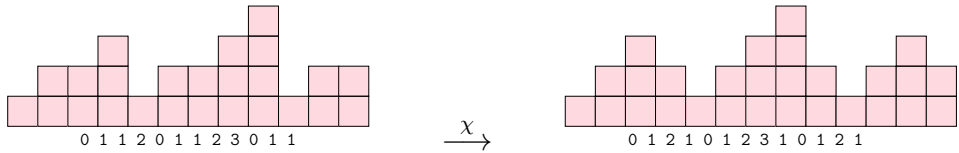
\begin{figure}[htb]
    \centering
   \scalebox{0.4}{\begin{tikzpicture}
% Définir la couleur de remplissage
%\definecolor{lightpink}{rgb}{1.0, 0.71, 0.76}
\definecolor{lightpink}{rgb}{1.0, 0.85, 0.88}
% Dessiner les cellules du polyomino avec leurs délimitations
% Colonne 1 (1 cellule)
\fill[lightpink] (0,0) rectangle (1,1);
\draw[black] (0,0) rectangle (1,1);
% Colonne 2 (2 cellules)
\fill[lightpink] (1,0) rectangle (2,1);
\fill[lightpink] (1,1) rectangle (2,2);
\draw[black] (1,0) rectangle (2,1);
\draw[black] (1,1) rectangle (2,2);
% Colonne 3 (3 cellules)
\fill[lightpink] (2,0) rectangle (3,1);
\fill[lightpink] (2,1) rectangle (3,2);
\draw[black] (2,0) rectangle (3,1);
\draw[black] (2,1) rectangle (3,2);
% Colonne 4 (1 cellule)
\fill[lightpink] (3,0) rectangle (4,1);
\fill[lightpink] (3,1) rectangle (4,2);
\fill[lightpink] (3,2) rectangle (4,3);
\draw[black] (3,0) rectangle (4,1);
\draw[black] (3,1) rectangle (4,2);
\draw[black] (3,2) rectangle (4,3);

% Colonne 5 (1 cellule)
\fill[lightpink] (4,0) rectangle (5,1);
\draw[black] (4,0) rectangle (5,1);

% Colonne 6 (1 cellule)
\fill[lightpink] (5,0) rectangle (6,1);
\fill[lightpink] (5,1) rectangle (6,2);
\draw[black] (5,0) rectangle (6,1);
\draw[black] (5,1) rectangle (6,2);

% Colonne 7 (1 cellule)
\fill[lightpink] (6,0) rectangle (7,1);
\fill[lightpink] (6,1) rectangle (7,2);
\draw[black] (6,0) rectangle (7,1);
\draw[black] (6,1) rectangle (7,2);

% Colonne 8 (1 cellule)
\fill[lightpink] (7,0) rectangle (8,1);
\fill[lightpink] (7,1) rectangle (8,2);
\fill[lightpink] (7,2) rectangle (8,3);
\draw[black] (7,0) rectangle (8,1);
\draw[black] (7,1) rectangle (8,2);
\draw[black] (7,2) rectangle (8,3);

% Colonne 9 (1 cellule)
\fill[lightpink] (8,0) rectangle (9,1);
\fill[lightpink] (8,1) rectangle (9,2);
\fill[lightpink] (8,2) rectangle (9,3);
\fill[lightpink] (8,3) rectangle (9,4);
\draw[black] (8,0) rectangle (9,1);
\draw[black] (8,1) rectangle (9,2);
\draw[black] (8,2) rectangle (9,3);
\draw[black] (8,3) rectangle (9,4);
% Colonne 10 (1 cellule)
\fill[lightpink] (9,0) rectangle (10,1);
\draw[black] (9,0) rectangle (10,1);

% Colonne 11 (1 cellule)
\fill[lightpink] (10,0) rectangle (11,1);
\fill[lightpink] (10,1) rectangle (11,2);
\draw[black] (10,0) rectangle (11,1);
\draw[black] (10,1) rectangle (11,2);

% Colonne 12 (1 cellule)
\fill[lightpink] (11,0) rectangle (12,1);
\fill[lightpink] (11,1) rectangle (12,2);
\draw[black] (11,0) rectangle (12,1);
\draw[black] (11,1) rectangle (12,2);

\node at (6,-0.45){\Large \texttt{\Large \texttt{0~1~1~2~0~1~1~2~3~0~1~1}}};
\end{tikzpicture}} \quad\quad
$\overset{\chi}{\longrightarrow}$
 \quad\quad  \scalebox{0.4}{\begin{tikzpicture}
% Définir la couleur de remplissage
%\definecolor{lightpink}{rgb}{1.0, 0.71, 0.76}
\definecolor{lightpink}{rgb}{1.0, 0.85, 0.88}
% Dessiner les cellules du polyomino avec leurs délimitations
% Colonne 1 (1 cellule)
\fill[lightpink] (0,0) rectangle (1,1);
\draw[black] (0,0) rectangle (1,1);
% Colonne 2 (2 cellules)
\fill[lightpink] (1,0) rectangle (2,1);
\fill[lightpink] (1,1) rectangle (2,2);
\draw[black] (1,0) rectangle (2,1);
\draw[black] (1,1) rectangle (2,2);
% Colonne 3 (3 cellules)
\fill[lightpink] (2,0) rectangle (3,1);
\fill[lightpink] (2,1) rectangle (3,2);\fill[lightpink] (2,2) rectangle (3,3);
\draw[black] (2,0) rectangle (3,1);
\draw[black] (2,1) rectangle (3,2);\draw[black] (2,2) rectangle (3,3);
% Colonne 4 (1 cellule)
\fill[lightpink] (3,0) rectangle (4,1);
\fill[lightpink] (3,1) rectangle (4,2);

\draw[black] (3,0) rectangle (4,1);
\draw[black] (3,1) rectangle (4,2);

% Colonne 5 (1 cellule)
\fill[lightpink] (4,0) rectangle (5,1);
\draw[black] (4,0) rectangle (5,1);

% Colonne 6 (1 cellule)
\fill[lightpink] (5,0) rectangle (6,1);
\fill[lightpink] (5,1) rectangle (6,2);
\draw[black] (5,0) rectangle (6,1);
\draw[black] (5,1) rectangle (6,2);

% Colonne 7 (1 cellule)
\fill[lightpink] (6,0) rectangle (7,1);
\fill[lightpink] (6,1) rectangle (7,2);\fill[lightpink] (6,2) rectangle (7,3);
\draw[black] (6,0) rectangle (7,1);
\draw[black] (6,1) rectangle (7,2);\draw[black] (6,2) rectangle (7,3);

% Colonne 8 (1 cellule)
\fill[lightpink] (7,0) rectangle (8,1);
\fill[lightpink] (7,1) rectangle (8,2);
\fill[lightpink] (7,2) rectangle (8,3);\fill[lightpink] (7,3) rectangle (8,4);
\draw[black] (7,0) rectangle (8,1);
\draw[black] (7,1) rectangle (8,2);
\draw[black] (7,2) rectangle (8,3);\draw[black] (7,3) rectangle (8,4);

% Colonne 9 (1 cellule)
\fill[lightpink] (8,0) rectangle (9,1);
\fill[lightpink] (8,1) rectangle (9,2);

\draw[black] (8,0) rectangle (9,1);
\draw[black] (8,1) rectangle (9,2);

% Colonne 10 (1 cellule)
\fill[lightpink] (9,0) rectangle (10,1);
\draw[black] (9,0) rectangle (10,1);

% Colonne 11 (1 cellule)
\fill[lightpink] (10,0) rectangle (11,1);
\fill[lightpink] (10,1) rectangle (11,2);
\draw[black] (10,0) rectangle (11,1);
\draw[black] (10,1) rectangle (11,2);

% Colonne 12 (1 cellule)
\fill[lightpink] (11,0) rectangle (12,1);
\fill[lightpink] (11,1) rectangle (12,2);\fill[lightpink] (11,2) rectangle (12,3);
\draw[black] (11,0) rectangle (12,1);
\draw[black] (11,1) rectangle (12,2);\draw[black] (11,2) rectangle (12,3);

% Colonne 12 (1 cellule)
\fill[lightpink] (12,0) rectangle (13,1);
\fill[lightpink] (12,1) rectangle (13,2);
\draw[black] (12,0) rectangle (13,1);
\draw[black] (12,1) rectangle (13,2);

\node at (6,-0.45){\Large\texttt{0~1~2~1~0~1~2~3~1~0~1~2~1}};
\end{tikzpicture}}

\caption{The bijection $\chi$ applied on $w=\texttt{011201123011}$. }\label{fig5}
\end{figure}

Notice that the area and the last symbol statistics do not transfer well through this bijection. For instance, when $w=\texttt{011201123011}$, we have 
$$
    \chi(\texttt{011201123011})=\texttt{0} (1+\chi(\texttt{00})) \chi(\texttt{01123011})=\texttt{0121}\texttt{0}(1+\chi(\texttt{001}))\chi(\texttt{011})=\texttt{0121}\texttt{0}\texttt{1231}\texttt{0121}
$$
and the semiperimeter of $w$ and $\chi(w)$ are 19 and 21, respectively (see Figure~\ref{fig5} for an illustration of this example).

Using relation \eqref{fq}, we deduce the following.
\begin{theorem}\label{smp}
The generating function for the number of nonempty $(\geq,\geq)$-Catalan words  according to the length, the semiperimeter, and the value of the last letter is given by
\[
C^{\geq}(x;p,1;v)=\frac{1+(1-2v)p^2x-2vp^3x^2-\sqrt{1-2p^2x+(p^4-4p^3)x^2}}{2((1-v)-p^2v(1-v)x+p^3v^2x^2)}.
\]
\end{theorem}
The first terms of the series expansion of $C^{\geq}(x;p,1;v)$ are
\begin{align*}
&p^2x+(p^3+p^4v)x^2+(p^5+2p^5v+p^6v^2)x^3+(p^6+p^6v+\mathbf{p^7+2p^7v+3p^7v^2}+p^8v^3)x^4\\
& \quad \quad \quad + (3p^8+p^9+(4p^8+2p^{9})v+3(p^8+p^9)v^2+4p^9v^3+p^{10}v^4)x^5+O(x^5).
\end{align*}

In Figure~\ref{fig4}, we can check that there are $1$, $2$, and $3$ $(\geq,\geq)$- polyominoes of length $4$ and semiperimeter $7$, whose value of the last symbol (ie. height of last column minus 1) is $0$, $1$, and $2$, respectively.

For $0\leq k\leq n-1$, let $c(n, k)$ denote the number of elements in the set $\mathcal{C}_{(\geq,\geq)}(n,k)$ of Catalan words avoiding the pattern $(\geq,\geq)$ and ending with the value $k$. Considering the two cases (i) and (ii) in the proof of Theorem~\ref{feq}, we can easily obtain the following.

\begin{corollary}For all $0 \leq k \leq n-1$, we have
\begin{align}
&c(n+1, 0)=c(n-1, 0)+c(n-1,1)+\cdots +c(n-1, n-2), \label{gc0}\\
&c(n+1, k+1)=c(n, k)+c(n-1, k)+\cdots +c(n-1, n-2),\label{gc}
\end{align}
where $c(1, 0)=c(2, 0)=c(2, 1)=1$.
\end{corollary}
The first few rows of the matrix $[c(n,k)]_{n\geq 1,k\geq0}$ are
\[
[c(n,k)]_{n\geq 1,k\geq0}=\left[
                            \begin{array}{cccccccc}
                              1 & 0 & 0 & 0 & 0  & 0 &0&\cdots \\
                              1 & 1 & 0 & 0 & 0 & 0 &0&\cdots\\
                              1 & 2 & 1 & 0 & 0 & 0 &0&\cdots\\
                              2 & 3 & 3 & 1 & 0 &0 &0&\cdots \\
                              4 & 6 & 6 & 4 & 1 & 0 &0&\cdots \\
                              9 & 13& 13& 10& 5 & 1 & 0&\cdots\\
                            21 & 30& 30& 24& 15 & 6 & 1&\cdots\\
                              \vdots & \vdots & \vdots & \vdots & \vdots &\vdots &\vdots &\ddots
                            \end{array}
                          \right].
\]
It is clear that, for $n\geq1$, $m_{n}=\sum_{k\geq0} c(n, k)$ and $m_n=c(n+2, 0)$, where $m_n$ is the $n$-th Motzkin number (see Introduction).

\medskip

Let $s(n,i)$ denote the sum of the semiperimeters of all $(\geq,\geq)$-polyominoes of length $n$ whose last column has height $i$. The first few values are
\[
[s(n,i)]_{n,i\geq 1}=\left[
                            \begin{array}{ccccccc}
                              2 & 0 & 0 & 0 & 0  &0&\cdots \\
                              3 & 4 & 0 & 0 & 0 &0&\cdots\\
                              5 & 10 & 6 & 0 & 0 &0&\cdots\\
                              13 & 20 & 21 & 8 & 0 &0&\cdots \\
                              33 & 50 & 51 & 36 & 10& 0&\cdots \\
                              89&130& 132&104&55&12&\cdots\\
                              \vdots & \vdots & \vdots & \vdots & \vdots & \vdots & \ddots
                            \end{array}
                          \right].
\]

\medskip  
Notice that from the decomposition of a $(\geq,\geq)$-Catalan word given in the proof of Theorem~\ref{feq} (see Figures~\ref{f14} and \ref{f24}), we have for $n\geq 3$ and $2\leq i \leq n-1$,
\[
s(n,i)=s(n-1,i-1)+2c(n-1,i-2)+\sum_{\substack{k\geq i-1 \\ k \neq i}}^{n-2}(s(n-2,k)+3c(n-2,k-1)).
\]

If we consider the difference $s(n,i)-s(n,i-1)$, then for $n\geq 3$ and $2\leq i \leq n-1$, we obtain the recurrence relation
\begin{align*}
& s(n,i)=s(n,i-1)+s(n-1,i-1)+s(n-2,i-1)-s(n-1,i-2) \\
&\quad \quad \quad -s(n-2,i)-s(n-2,i-2)+2(c(n-1,i-2)-c(n-1,i-3)\\
&\quad \quad \quad +3(c(n-2,i-2)-c(n-2,i-1)-c(n-2,i-2)).
\end{align*}

Let $s(n)$ be the total sum of semiperimeters over all $(\geq,\geq)$-polyominoes of length $n$. By Theorem \ref{smp} with $v=1$, and after differentiating at $p=1$, we deduce:
\begin{corollary}
The generating function of the sequence $s(n)$ is given by
\[
\frac{3-4x-5x^2+(x-3)\sqrt{1-2x-3x^2}}{2x^2\sqrt{1-2x-3x^2}}.
\]
An asymptotic approximation of $s(n)$ is given by

$$\frac{5\sqrt{3}}{2\sqrt{\pi}}\cdot \frac{3^n}{\sqrt{n}}.$$

\label{cor34}
\end{corollary}

The first few values for $1\leq n \leq 10$
\[
 2,\quad 7,\quad 21,\quad 62,\quad 180,\quad 522,\quad 1512,\quad 4384,\quad 12726,\quad 36995.
\]
This sequence does not appear in the OEIS. The expected value of the semiperimeter is $5n/3$. In Corollary~\ref{corsn} we give an explicit relation for the sequence $s(n)$ using the central trinomial coefficient $T_n$ (see Section \ref{sec:2}).

\begin{corollary}\label{corsn}The total sum of  semiperimeters over all $(\geq,\geq)$-polyominoes of length $n$ is given by
$$s(n)=\frac{1}{2}\cdot\left(-5 T_n -4 T_{n+1}+3T_{n+2}\right).$$
\end{corollary}
\begin{proof} We decompose the generating function of Corollary~\ref{cor34} as $$\frac{3-4x-5x^2+(x-3)\sqrt{1-2x-3x^2}}{2x^2\sqrt{1-2x-3x^2}}={\frac {x-3}{2{x}^{2}}}-\frac {5\,{x}^{2}+4\,x-3}{2x^2}\cdot T(x).
$$
Using relation \eqref{gft}, and comparing the $n$-th coefficient, we obtain the desired result.
\end{proof}

\section{The area statistic}
\label{sec:4}

The goal of this section is to analyze the area statistic on  $(\geq,\geq)$-polyominoes. Considering Theorem~\ref{feq} with $p=1$, we obtain
\begin{align}\label{f1}
C^{\geq}(x;1,q;v)&=qx+q^2x^2+\frac{q^3x^2}{1-qv}C^{\geq}(x;1,q;q)+q^2xvC^{\geq}(x;1,q;qv)\\ \notag
&\quad \quad\quad \quad\quad \quad \quad\quad\quad\quad\quad \quad\quad \quad\quad\quad \quad  -\frac{q^5 x^2v^2}{1-qv}C^{\geq}(x;1,q;q^2v).
\end{align}
Using this relation, we do not succeed to obtain a closed form for the generating function  $C^{\geq}(x;1,q;v)$. So, we will proceed differently from the previous section.

Let $\mathcal{B}$ be the set consisting of $(\geq,\geq)$-Catalan words $w\in \mathcal{C}$ such that $w$ does not end with $ab$ such that $a\geq b$.  Then we have the following.

\begin{theorem}\label{thb} The generating function $B(x,q)$ according to the length and the area for the number of nonempty Catalan words $w$ avoiding $(\geq,\geq)$  such that $w$ does not end with $ab$, $a\geq b$ is given by 
\begin{align*}
B(x,q)=\frac{\sum_{j\geq1}(-1)^{j-1}q^{j}x^j\prod_{i=1}^{j-1}\frac{1-q^i+q^{2i}}{1-q^i}}
{1-\sum_{j\geq1}(-1)^{j-1}\frac{q^{j}x^j}{1-q^j}\prod_{i=1}^{j-1}\frac{1-q^i+q^{2i}}{1-q^i}}.
\end{align*}
In term of infinite continued fraction, we have
\begin{align*}
1+B(x,q)&=\cfrac{1}{1-\cfrac{qx}{1+qx-\cfrac{(1+qx)q^2x}{1+q^2x-\cfrac{(1+q^2x)q^3x}{\ddots}}}}.
\end{align*}
\end{theorem}
\begin{proof} Let us define recursively  the bijection $\psi$ from $\mathcal{B}$ to the set of ($\neq$)-Catalan words as follows: 
$$\psi(w)=\left\{\begin{array}{ll}
\epsilon & \mbox{ if } w=\epsilon, \\
\texttt{0}(1+\psi(u)) \psi(v)& \mbox{ if } w=\texttt{0}(1+u)v \mbox{ and } u\neq\epsilon,\\
\psi(u)\texttt{0}& \mbox{ if } w=\texttt{0}u.\\
\end{array}\right.
$$
It is clear that $\psi$ is a bijection that transports the length and the area statistics. Indeed, if $w=\texttt{0}(1+u)v$, $v\neq\epsilon$,  then we have 
$\area(\psi(\texttt{0}(1+u)v))=\area(\texttt{0}(1+\psi(u)))+\area(\psi(v))$, and by induction and using the fact that $\psi$ preserves the length, we have  $\area(\texttt{0}(1+\psi(u)))+\area(\psi(v))=\area(\texttt{0}(1+u))+\area(v)=\area(w)$ .
Considering this bijection, we directly obtain the result using Theorem~5.1 in \cite{BKR}.
\end{proof}
The first terms of the series expansion of $B(x,q)$ are 
\begin{align*} q x + q^3 x^2 + (q^4 + q^6) x^3 &+ (q^6 + q^7 + q^8 +     q^{10}) x^4+\\
    &+ (q^7 + 3 q^9 + 2 q^{11} + q^{12} + q^{13} + q^{15}) x^5 + O(x^6).
\end{align*}

\begin{theorem}\label{thmaire}The generating function for the number of nonempty $(\geq,\geq)$-Catalan words according to the length and area is given by $$C^{\geq}(x;1,q;1)= \sum_{i=1}^{\infty} x^i q^{\frac{i(i+1)}{2}} \prod_{j=0}^{i-1} \left(1+\mathit{B} \left(x q^j, q \right)\right).
$$
 \end{theorem}
 \begin{proof} Any nonempty Catalan word avoiding $(\geq,\geq)$ is of the form either $u0$ or $u0(1+v)$ where $u\in\mathcal{B}$ and $v$ is a $(\geq,\geq)$-Catalan word, which generate the functional equation  
$$C^{\geq}(x;1,q;1)=xq(1+B(x,q))+xqC^{\geq}(xq;1,q;1)(1+B(x,q)).$$
Iterating this equation we obtain the desired result.
 \end{proof}
The first terms of the series expansion of $C^{\geq}(x;1,q;1)$ are 

\begin{align*}qx + (q^3 + q^2)x^2 + (q^6 + q^5 + 2q^4)x^3 + (q^{10} + q^9 + 2q^8 + 2q^7 + 2q^6 + q^5)x^4 +\quad\quad\quad\quad\quad\quad\\ (q^{15} + q^{14} + 2q^{13} + 2q^{12} + 4q^{11} + 2q^{10} + \mathbf{5q^9} + 2q^8 + 2q^7)x^5 +\quad\quad\quad\quad\\ (q^{21} + q^{20} + 2q^{19} + 2q^{18} + 4q^{17} + 4q^{16} + 6q^{15} + 5q^{14} + 7q^{13} + 6q^{12} + 5q^{11} + 5q^{10} + 2q^9 + q^8)x^6 +O(x^7)
\end{align*}

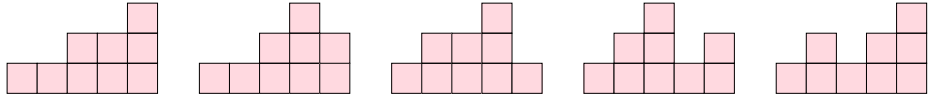
\begin{figure}[htb]
\centering
\scalebox{0.4}{\begin{tikzpicture}
% Définir la couleur de remplissage
%\definecolor{lightpink}{rgb}{1.0, 0.71, 0.76}
\definecolor{lightpink}{rgb}{1.0, 0.85, 0.88}
% Dessiner les cellules du polyomino avec leurs délimitations
% Colonne 1 (1 cellule)
\fill[lightpink] (0,0) rectangle (1,1);
\draw[black] (0,0) rectangle (1,1);
% Colonne 2 (2 cellules)
\fill[lightpink] (1,0) rectangle (2,1);
%\fill[lightpink] (1,1) rectangle (2,2);
\draw[black] (1,0) rectangle (2,1);
%\draw[black] (1,1) rectangle (2,2);
% Colonne 3 (3 cellules)
\fill[lightpink] (2,0) rectangle (3,1);
\fill[lightpink] (2,1) rectangle (3,2);
%\fill[lightpink] (2,2) rectangle (3,3);
\draw[black] (2,0) rectangle (3,1);
\draw[black] (2,1) rectangle (3,2);
%\draw[black] (2,2) rectangle (3,3);
% Colonne 4 (1 cellule)
\fill[lightpink] (3,0) rectangle (4,1);
\fill[lightpink] (3,1) rectangle (4,2);
\draw[black] (3,0) rectangle (4,1);
\draw[black] (3,1) rectangle (4,2);
% Colonne 5 (1 cellule)
\fill[lightpink] (4,0) rectangle (5,1);
\fill[lightpink] (4,1) rectangle (5,2);
\fill[lightpink] (4,2) rectangle (5,3);
\draw[black] (4,0) rectangle (5,1);
\draw[black] (4,1) rectangle (5,2);
\draw[black] (4,2) rectangle (5,3);
\end{tikzpicture}}
\quad 
\scalebox{0.4}{\begin{tikzpicture}
% Définir la couleur de remplissage
%\definecolor{lightpink}{rgb}{1.0, 0.71, 0.76}
\definecolor{lightpink}{rgb}{1.0, 0.85, 0.88}
% Dessiner les cellules du polyomino avec leurs délimitations
% Colonne 1 (1 cellule)
\fill[lightpink] (0,0) rectangle (1,1);
\draw[black] (0,0) rectangle (1,1);
% Colonne 2 (2 cellules)
\fill[lightpink] (1,0) rectangle (2,1);
%\fill[lightpink] (1,1) rectangle (2,2);
\draw[black] (1,0) rectangle (2,1);
%\draw[black] (1,1) rectangle (2,2);
% Colonne 3 (3 cellules)
\fill[lightpink] (2,0) rectangle (3,1);
\fill[lightpink] (2,1) rectangle (3,2);
%\fill[lightpink] (2,2) rectangle (3,3);
\draw[black] (2,0) rectangle (3,1);
\draw[black] (2,1) rectangle (3,2);
%\draw[black] (2,2) rectangle (3,3);
% Colonne 4 (1 cellule)
\fill[lightpink] (3,0) rectangle (4,1);
\fill[lightpink] (3,1) rectangle (4,2);
\fill[lightpink] (3,2) rectangle (4,3);
\draw[black] (3,0) rectangle (4,1);
\draw[black] (3,1) rectangle (4,2);
\draw[black] (3,2) rectangle (4,3);
% Colonne 4 (1 cellule)
\fill[lightpink] (4,0) rectangle (5,1);
\draw[black] (4,0) rectangle (5,1);
\fill[lightpink] (4,1) rectangle (5,2);
\draw[black] (4,1) rectangle (5,2);
\end{tikzpicture}}
\quad 
\scalebox{0.4}{\begin{tikzpicture}
% Définir la couleur de remplissage
%\definecolor{lightpink}{rgb}{1.0, 0.71, 0.76}
\definecolor{lightpink}{rgb}{1.0, 0.85, 0.88}
% Dessiner les cellules du polyomino avec leurs délimitations
% Colonne 1 (1 cellule)
\fill[lightpink] (0,0) rectangle (1,1);
\draw[black] (0,0) rectangle (1,1);
% Colonne 2 (2 cellules)
\fill[lightpink] (1,0) rectangle (2,1);
\fill[lightpink] (1,1) rectangle (2,2);
\draw[black] (1,0) rectangle (2,1);
\draw[black] (1,1) rectangle (2,2);
% Colonne 3 (3 cellules)
\fill[lightpink] (2,0) rectangle (3,1);
\fill[lightpink] (2,1) rectangle (3,2);
%\fill[lightpink] (2,2) rectangle (3,3);
\draw[black] (2,0) rectangle (3,1);
\draw[black] (2,1) rectangle (3,2);
%\draw[black] (2,2) rectangle (3,3);
% Colonne 4 (1 cellule)
\fill[lightpink] (3,0) rectangle (4,1);
\fill[lightpink] (3,1) rectangle (4,2);
\fill[lightpink] (3,2) rectangle (4,3);
\draw[black] (3,0) rectangle (4,1);
\draw[black] (3,1) rectangle (4,2);
\draw[black] (3,2) rectangle (4,3);
% Colonne 4 (1 cellule)
\fill[lightpink] (4,0) rectangle (5,1);
\draw[black] (4,0) rectangle (5,1);
\end{tikzpicture}}
\quad \scalebox{0.4}{\begin{tikzpicture}
% Définir la couleur de remplissage
%\definecolor{lightpink}{rgb}{1.0, 0.71, 0.76}
\definecolor{lightpink}{rgb}{1.0, 0.85, 0.88}
% Dessiner les cellules du polyomino avec leurs délimitations
% Colonne 1 (1 cellule)
\fill[lightpink] (0,0) rectangle (1,1);
\draw[black] (0,0) rectangle (1,1);
% Colonne 2 (2 cellules)
\fill[lightpink] (1,0) rectangle (2,1);
\fill[lightpink] (1,1) rectangle (2,2);
\draw[black] (1,0) rectangle (2,1);
\draw[black] (1,1) rectangle (2,2);
% Colonne 3 (3 cellules)
\fill[lightpink] (2,0) rectangle (3,1);
\fill[lightpink] (2,1) rectangle (3,2);
\fill[lightpink] (2,2) rectangle (3,3);
\draw[black] (2,0) rectangle (3,1);
\draw[black] (2,1) rectangle (3,2);
\draw[black] (2,2) rectangle (3,3);
% Colonne 4 (1 cellule)
\fill[lightpink] (3,0) rectangle (4,1);
%\fill[lightpink] (3,2) rectangle (4,3);
\draw[black] (3,0) rectangle (4,1);
%\draw[black] (3,2) rectangle (4,3);
% Colonne 5 (1 cellule)
\fill[lightpink] (4,0) rectangle (5,1);
\fill[lightpink] (4,1) rectangle (5,2);
\draw[black] (4,0) rectangle (5,1);
\draw[black] (4,1) rectangle (5,2);
\end{tikzpicture}}
\quad \scalebox{0.4}{\begin{tikzpicture}
% Définir la couleur de remplissage
%\definecolor{lightpink}{rgb}{1.0, 0.71, 0.76}
\definecolor{lightpink}{rgb}{1.0, 0.85, 0.88}
% Dessiner les cellules du polyomino avec leurs délimitations
% Colonne 1 (1 cellule)
\fill[lightpink] (0,0) rectangle (1,1);
\draw[black] (0,0) rectangle (1,1);
% Colonne 2 (2 cellules)
\fill[lightpink] (1,0) rectangle (2,1);
\fill[lightpink] (1,1) rectangle (2,2);
\draw[black] (1,0) rectangle (2,1);
\draw[black] (1,1) rectangle (2,2);
% Colonne 3 (3 cellules)
\fill[lightpink] (2,0) rectangle (3,1);

\draw[black] (2,0) rectangle (3,1);

% Colonne 4 (1 cellule)
\fill[lightpink] (3,0) rectangle (4,1);
\fill[lightpink] (3,1) rectangle (4,2);

\draw[black] (3,0) rectangle (4,1);
\draw[black] (3,1) rectangle (4,2);
% Colonne 5 (1 cellule)
\fill[lightpink] (4,0) rectangle (5,1);
\fill[lightpink] (4,1) rectangle (5,2);
\fill[lightpink] (4,2) rectangle (5,3);
\draw[black] (4,0) rectangle (5,1);
\draw[black] (4,1) rectangle (5,2);
\draw[black] (4,2) rectangle (5,3);
\end{tikzpicture}}
\caption{The $5$ $(\geq,\geq)$-polyominoes of area 9 and length $5$. }\label{figaire}
\end{figure}

Let $u(n, i)$ denote the total area of the $(\geq,\geq)$-polyominoes of length $n$ that end with a column of height $i$. The first
few values are
\[
[u(n,i)]_{n,i\geq 1}=\left[
                            \begin{array}{cccccc}
                              1 & 0 & 0 & 0 & 0  &\cdots \\
                              2 & 3 & 0 & 0 & 0 &\cdots\\
                              4 & 9 & 6 & 0 & 0 &\cdots\\
                              12 & 20 & 24 & 10 & 0 &\cdots \\
                              35& 55 & 63 & 50 & 15&\cdots\\
                              \vdots & \vdots & \vdots & \vdots & \vdots & \ddots
                            \end{array}
                          \right].
\]

From the decomposition given in Figures~\ref{f14} and \ref{f24}, we have for $n\geq2$ and $2 \leq i \leq n$,
\[
u(n,i)=u(n-1,i-1)+(i+1)c(n-1,i-2)+\sum_{k=i}^{n-2}(u(n-2,k-1)+(i+k+2)c(n-2,k-2)).
\]
If we consider the difference $u(n,i)-u(n,i-1)$, then for $n\geq2$ and $3\leq i\leq n$, we obtain
the recurrence relation

\begin{align*}
&u(n,i)=u(n,i-1)+u(n-1,i-1)-u(n-1,i-2)+ u(n-2,n-3)\\
&\quad \quad \quad -u(n-2,i-2)+(n+i)c(n-2,n-4)-(2i+1)c(n-2,i-3)\\
&\quad \quad \quad +(i+1)c(n-1,i-2)-i.c(n-1,i-3)+\sum_{k=i-1}^{n-2}c(n-2,k-2).
\end{align*}

Let $u(n)$ be the total sum of area over all $(\geq,\geq)$-polyominoes of length $n$. The first few values for $1\leq n\leq 10$ are
\[
1,\quad 5,\quad 19,\quad 66,\quad 218,\quad 701,\quad 2215,\quad 6919,\quad 21438,\quad 66034.
\]
The
following theorem gives a combinatorial formula to calculate the sequence $u(n)$.
\begin{theorem} The generating function for the total sum of area over all $(\geq,\geq)$-Catalan words with respect to the length is 
\begin{align*}
U(x,1)=\frac{\left(-2 x^{3}+x^{2}+3 x -1\right) \sqrt{-3 x^{2}-2 x +1}+2 x^{4}+7 x^{3}-4 x +1}{6 x^{6}+4 x^{5}-2 x^{4}},\end{align*}
and the coefficient of $x^n$ is given given by
\[
u(n)=\frac{1}{2}\cdot\left(3^{n+1}+2T_{n+1}-T_{n+2}-3T_{n+3}+T_{n+4}\right).
\]
An asymptotic approximation is $3^{n+1}/2.$
\end{theorem}
\begin{proof}
Let $U(x;v)=\left.\frac{\partial}{\partial q}C^{\geq}(x;1,q;v) \right|_{q=1}$. Then by differentiating \eqref{f1} with respect to $q$, we obtain
\begin{align*}
U(x;v)&= x+2x^2+\frac{x^2(3-2v)(M(x)-1)}{(1-v)^2}+\frac{x^2}{1-v}\left(U(x;1)+\frac{\partial}{\partial v}C^{\geq}(x;1,1;v)|_{v=1}\right)\\
& \quad \quad +\left(2xv-\frac{x^2v^2(5-4v)}{(1-v)^2}\right)C^{\geq}(x;1,1;v)+xv\left(U(x;v)+v\frac{\partial}{\partial v}C^{\geq}(x;1,1;v)\right)\\
&\quad \quad -\frac{x^2v^2}{1-v}\left(U(x;v)+2v\frac{\partial}{\partial v}C^{\geq}(x;1,1;v)\right).
\end{align*}

Notice that with relation \eqref{gf3}, we obtain
\begin{align*}
&\frac{\partial}{\partial v} C^{\geq}(x;1,1;v)=\frac{\partial}{\partial v}\left(\frac{x(1-v)-x^2v+x^2M(x)}{1-v-xv(1-v)+x^2v^2}\right) \label{f1}\\ 
&\quad =\frac{\left(x +1\right) \left(\left(x v -\frac{1}{2}\right) \sqrt{-3 x^{2}-2 x +1}+\frac{1}{2}+v^{2} x^{3}+\left(v^{2}-v \right) x^{2}+\left(-v -\frac{1}{2}\right) x \right)}{\left(v^{2} x^{2}+\left(v^{2}-v \right) x -v +1\right)^{2}}.
\end{align*}

Now, we group the terms $U(x;v)$ and we rewrite the equation so that the left-hand side is 
$$\frac{U \! \left(x , v\right) \left(1-v -x v \left(-v +1\right)+v^{2} x^{2}\right)^3}{1-v}.$$

By twice differentiating this equation with respect to $v$ and canceling the factor $1-v-xv(1-v)+x^2v^2$ with $$v={\frac {1+x-\sqrt {-3\,{x}^{2}-2\,x+1}}{2x \left( x+1 \right) }}=\frac{1+xM(x)}{1+x},$$ we obtain that
\begin{align*}
U(x,1)&=\frac{\left(-2 x^{3}+x^{2}+3 x -1\right) \sqrt{-3 x^{2}-2 x +1}+2 x^{4}+7 x^{3}-4 x +1}{6 x^{6}+4 x^{5}-2 x^{4}}
. \end{align*}
The following decomposition
\begin{align*}
U(x,1)&=\frac{2 x^{3}+5 x^{2}-5 x +1}{2x^4(3x-1)} +\frac{\left(2 x^{3}-x^{2}-3 x +1\right) }{2 x^{4}}T \! \left(x \right)\end{align*}
allows to obtain the desired result (see Section \ref{sec:2} for the definition of $T(x)$).
\end{proof}
%\subsection{Link with ?? coefficients}
%Let $h(n,i)$ be the total number of cells of height $i$ in $\mathcal{C}_{(\geq, %\geq)}(n)$. The first
%few values are
%\[
%[h(n,i)]_{n,i\geq 1}=\left[
%                            \begin{array}{cccccc}
%                              1 & 0 & 0 & 0   &\cdots \\
%                              1 & 4 & 0 & 0  &\cdots\\
%                              1 & 6 & 12 & 0  &\cdots\\
%                              1 & 8 & 21 & 36  &\cdots \\
%                              \vdots & \vdots & \vdots & \vdots  & \ddots
%                            \end{array}
%                          \right].
%\]

%\medskip

%Let $w$ be a  Catalan word that avoid the pattern $(\geq, \geq)$. We denote by %$h_i(w)$ the number of cells of height $i$ in the Motzkin polyomino associated %with $w$. We introduce the following generating functions
%\[
%H_i(x, q) := 1+\sum_{w\in \mathcal{C}_{(\geq, \geq)}(n)}x^{|w|}q^{h_i(w)}
%\]
%and
%\[
%B_i(x)=\left. \frac{\partial H_i(x, q) }{\partial q}\right|_{q=1}.
%\]
%From the definition it is clear that $[x^n]B_i(x)=h(n,i)$.

\section{The interior points statistic}
\label{sec:5}

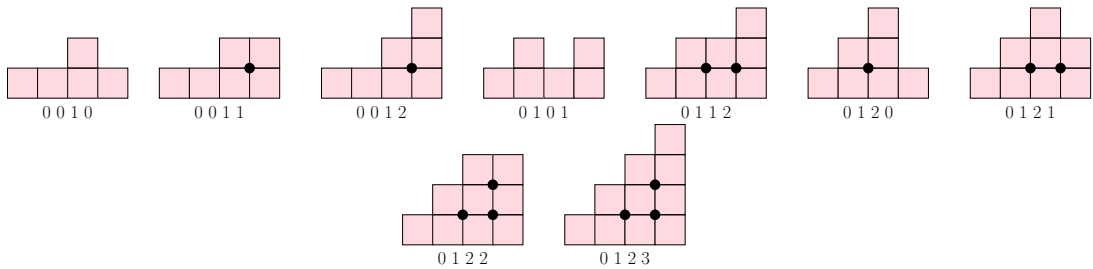
\begin{figure}[htb]
\centering
\scalebox{0.4}{\begin{tikzpicture}
% Définir la couleur de remplissage
%\definecolor{lightpink}{rgb}{1.0, 0.71, 0.76}
\definecolor{lightpink}{rgb}{1.0, 0.85, 0.88}
% Dessiner les cellules du polyomino avec leurs délimitations
% Colonne 1 (1 cellule)
\fill[lightpink] (0,0) rectangle (1,1);
\draw[black] (0,0) rectangle (1,1);
% Colonne 2 (2 cellules)
\fill[lightpink] (1,0) rectangle (2,1);
%\fill[lightpink] (1,1) rectangle (2,2);
\draw[black] (1,0) rectangle (2,1);
%\draw[black] (1,1) rectangle (2,2);
% Colonne 3 (3 cellules)
\fill[lightpink] (2,0) rectangle (3,1);
\fill[lightpink] (2,1) rectangle (3,2);
%\fill[lightpink] (2,2) rectangle (3,3);
\draw[black] (2,0) rectangle (3,1);
\draw[black] (2,1) rectangle (3,2);
%\draw[black] (2,2) rectangle (3,3);
% Colonne 4 (1 cellule)
\fill[lightpink] (3,0) rectangle (4,1);
\draw[black] (3,0) rectangle (4,1);
\node at (2,-0.45){\Large 0~0~1~0};
\end{tikzpicture}}%*****************
\quad 
\scalebox{0.4}{\begin{tikzpicture}
% Définir la couleur de remplissage
%\definecolor{lightpink}{rgb}{1.0, 0.71, 0.76}
\definecolor{lightpink}{rgb}{1.0, 0.85, 0.88}
% Dessiner les cellules du polyomino avec leurs délimitations
% Colonne 1 (1 cellule)
\fill[lightpink] (0,0) rectangle (1,1);
\draw[black] (0,0) rectangle (1,1);
% Colonne 2 (2 cellules)
\fill[lightpink] (1,0) rectangle (2,1);
%\fill[lightpink] (1,1) rectangle (2,2);
\draw[black] (1,0) rectangle (2,1);
%\draw[black] (1,1) rectangle (2,2);
% Colonne 3 (3 cellules)
\fill[lightpink] (2,0) rectangle (3,1);
\fill[lightpink] (2,1) rectangle (3,2);
%\fill[lightpink] (2,2) rectangle (3,3);
\draw[black] (2,0) rectangle (3,1);
\draw[black] (2,1) rectangle (3,2);
%\draw[black] (2,2) rectangle (3,3);
% Colonne 4 (1 cellule)
\fill[lightpink] (3,0) rectangle (4,1);
\fill[lightpink] (3,1) rectangle (4,2);
\draw[black] (3,0) rectangle (4,1);
\draw[black] (3,1) rectangle (4,2);\fill[black] (3,1) circle (5pt);
\node at (2,-0.45){\Large 0~0~1~1};\end{tikzpicture}}
\quad 
\scalebox{0.4}{\begin{tikzpicture}
% Définir la couleur de remplissage
%\definecolor{lightpink}{rgb}{1.0, 0.71, 0.76}
\definecolor{lightpink}{rgb}{1.0, 0.85, 0.88}
% Dessiner les cellules du polyomino avec leurs délimitations
% Colonne 1 (1 cellule)
\fill[lightpink] (0,0) rectangle (1,1);
\draw[black] (0,0) rectangle (1,1);
% Colonne 2 (2 cellules)
\fill[lightpink] (1,0) rectangle (2,1);
%\fill[lightpink] (1,1) rectangle (2,2);
\draw[black] (1,0) rectangle (2,1);
%\draw[black] (1,1) rectangle (2,2);
% Colonne 3 (3 cellules)
\fill[lightpink] (2,0) rectangle (3,1);
\fill[lightpink] (2,1) rectangle (3,2);
%\fill[lightpink] (2,2) rectangle (3,3);
\draw[black] (2,0) rectangle (3,1);
\draw[black] (2,1) rectangle (3,2);
%\draw[black] (2,2) rectangle (3,3);
% Colonne 4 (1 cellule)
\fill[lightpink] (3,0) rectangle (4,1);
\fill[lightpink] (3,1) rectangle (4,2);
\fill[lightpink] (3,2) rectangle (4,3);
\draw[black] (3,0) rectangle (4,1);
\draw[black] (3,1) rectangle (4,2);
\draw[black] (3,2) rectangle (4,3);
\node at (2,-0.45){\Large 0~0~1~2};\fill[black] (3,1) circle (5pt);
\end{tikzpicture}}
\quad 
\scalebox{0.4}{\begin{tikzpicture}
% Définir la couleur de remplissage
%\definecolor{lightpink}{rgb}{1.0, 0.71, 0.76}
\definecolor{lightpink}{rgb}{1.0, 0.85, 0.88}
% Dessiner les cellules du polyomino avec leurs délimitations
% Colonne 1 (1 cellule)
\fill[lightpink] (0,0) rectangle (1,1);
\draw[black] (0,0) rectangle (1,1);
% Colonne 2 (2 cellules)
\fill[lightpink] (1,0) rectangle (2,1);
\fill[lightpink] (1,1) rectangle (2,2);
\draw[black] (1,0) rectangle (2,1);
\draw[black] (1,1) rectangle (2,2);
% Colonne 3 (3 cellules)
\fill[lightpink] (2,0) rectangle (3,1);
%\fill[lightpink] (2,1) rectangle (3,2);
%\fill[lightpink] (2,2) rectangle (3,3);
\draw[black] (2,0) rectangle (3,1);
%\draw[black] (2,1) rectangle (3,2);
%\draw[black] (2,2) rectangle (3,3);
% Colonne 4 (1 cellule)
\fill[lightpink] (3,0) rectangle (4,1);
\fill[lightpink] (3,1) rectangle (4,2);
\draw[black] (3,0) rectangle (4,1);
\draw[black] (3,1) rectangle (4,2);
\node at (2,-0.45){\Large 0~1~0~1};
\end{tikzpicture}}
\quad 
\scalebox{0.4}{\begin{tikzpicture}
% Définir la couleur de remplissage
%\definecolor{lightpink}{rgb}{1.0, 0.71, 0.76}
\definecolor{lightpink}{rgb}{1.0, 0.85, 0.88}
% Dessiner les cellules du polyomino avec leurs délimitations
% Colonne 1 (1 cellule)
\fill[lightpink] (0,0) rectangle (1,1);
\draw[black] (0,0) rectangle (1,1);
% Colonne 2 (2 cellules)
\fill[lightpink] (1,0) rectangle (2,1);
\fill[lightpink] (1,1) rectangle (2,2);
\draw[black] (1,0) rectangle (2,1);
\draw[black] (1,1) rectangle (2,2);
% Colonne 3 (3 cellules)
\fill[lightpink] (2,0) rectangle (3,1);
\fill[lightpink] (2,1) rectangle (3,2);
%\fill[lightpink] (2,2) rectangle (3,3);
\draw[black] (2,0) rectangle (3,1);
\draw[black] (2,1) rectangle (3,2);
%\draw[black] (2,2) rectangle (3,3);
% Colonne 4 (1 cellule)
\fill[lightpink] (3,0) rectangle (4,1);
\fill[lightpink] (3,1) rectangle (4,2);
\fill[lightpink] (3,2) rectangle (4,3);
\draw[black] (3,0) rectangle (4,1);
\draw[black] (3,1) rectangle (4,2);
\draw[black] (3,2) rectangle (4,3);
\node at (2,-0.45){\Large 0~1~1~2};
\fill[black] (2,1) circle (5pt);\fill[black] (3,1) circle (5pt);
\end{tikzpicture}}
\quad 
\scalebox{0.4}{\begin{tikzpicture}
% Définir la couleur de remplissage
%\definecolor{lightpink}{rgb}{1.0, 0.71, 0.76}
\definecolor{lightpink}{rgb}{1.0, 0.85, 0.88}
% Dessiner les cellules du polyomino avec leurs délimitations
% Colonne 1 (1 cellule)
\fill[lightpink] (0,0) rectangle (1,1);
\draw[black] (0,0) rectangle (1,1);
% Colonne 2 (2 cellules)
\fill[lightpink] (1,0) rectangle (2,1);
\fill[lightpink] (1,1) rectangle (2,2);
\draw[black] (1,0) rectangle (2,1);
\draw[black] (1,1) rectangle (2,2);
% Colonne 3 (3 cellules)
\fill[lightpink] (2,0) rectangle (3,1);
\fill[lightpink] (2,1) rectangle (3,2);
\fill[lightpink] (2,2) rectangle (3,3);
\draw[black] (2,0) rectangle (3,1);
\draw[black] (2,1) rectangle (3,2);
\draw[black] (2,2) rectangle (3,3);
% Colonne 4 (1 cellule)
\fill[lightpink] (3,0) rectangle (4,1);
%\fill[lightpink] (3,1) rectangle (4,2);
%\fill[lightpink] (3,2) rectangle (4,3);
\draw[black] (3,0) rectangle (4,1);
%\draw[black] (3,1) rectangle (4,2);
%\draw[black] (3,2) rectangle (4,3);
\node at (2,-0.45){\Large 0~1~2~0};\fill[black] (2,1) circle (5pt);
\end{tikzpicture}}
\quad \scalebox{0.4}{\begin{tikzpicture}
% Définir la couleur de remplissage
%\definecolor{lightpink}{rgb}{1.0, 0.71, 0.76}
\definecolor{lightpink}{rgb}{1.0, 0.85, 0.88}
% Dessiner les cellules du polyomino avec leurs délimitations
% Colonne 1 (1 cellule)
\fill[lightpink] (0,0) rectangle (1,1);
\draw[black] (0,0) rectangle (1,1);
% Colonne 2 (2 cellules)
\fill[lightpink] (1,0) rectangle (2,1);
\fill[lightpink] (1,1) rectangle (2,2);
\draw[black] (1,0) rectangle (2,1);
\draw[black] (1,1) rectangle (2,2);
% Colonne 3 (3 cellules)
\fill[lightpink] (2,0) rectangle (3,1);
\fill[lightpink] (2,1) rectangle (3,2);
\fill[lightpink] (2,2) rectangle (3,3);
\draw[black] (2,0) rectangle (3,1);
\draw[black] (2,1) rectangle (3,2);
\draw[black] (2,2) rectangle (3,3);
% Colonne 4 (1 cellule)
\fill[lightpink] (3,0) rectangle (4,1);
\fill[lightpink] (3,1) rectangle (4,2);
%\fill[lightpink] (3,2) rectangle (4,3);
\draw[black] (3,0) rectangle (4,1);
\draw[black] (3,1) rectangle (4,2);
%\draw[black] (3,2) rectangle (4,3);
\node at (2,-0.45){\Large 0~1~2~1};\fill[black] (2,1) circle (5pt);\fill[black] (3,1) circle (5pt);
\end{tikzpicture}}
\quad \scalebox{0.4}{\begin{tikzpicture}
% Définir la couleur de remplissage
%\definecolor{lightpink}{rgb}{1.0, 0.71, 0.76}
\definecolor{lightpink}{rgb}{1.0, 0.85, 0.88}
% Dessiner les cellules du polyomino avec leurs délimitations
% Colonne 1 (1 cellule)
\fill[lightpink] (0,0) rectangle (1,1);
\draw[black] (0,0) rectangle (1,1);
% Colonne 2 (2 cellules)
\fill[lightpink] (1,0) rectangle (2,1);
\fill[lightpink] (1,1) rectangle (2,2);
\draw[black] (1,0) rectangle (2,1);
\draw[black] (1,1) rectangle (2,2);
% Colonne 3 (3 cellules)
\fill[lightpink] (2,0) rectangle (3,1);
\fill[lightpink] (2,1) rectangle (3,2);
\fill[lightpink] (2,2) rectangle (3,3);
\draw[black] (2,0) rectangle (3,1);
\draw[black] (2,1) rectangle (3,2);
\draw[black] (2,2) rectangle (3,3);
% Colonne 4 (1 cellule)
\fill[lightpink] (3,0) rectangle (4,1);
\fill[lightpink] (3,1) rectangle (4,2);
\fill[lightpink] (3,2) rectangle (4,3);
\draw[black] (3,0) rectangle (4,1);
\draw[black] (3,1) rectangle (4,2);
\draw[black] (3,2) rectangle (4,3);
\node at (2,-0.45){\Large 0~1~2~2};
\fill[black] (2,1) circle (5pt);
\fill[black] (3,1) circle (5pt);\fill[black] (3,2) circle (5pt);
\end{tikzpicture}}
\quad \scalebox{0.4}{\begin{tikzpicture}
% Définir la couleur de remplissage
%\definecolor{lightpink}{rgb}{1.0, 0.71, 0.76}
\definecolor{lightpink}{rgb}{1.0, 0.85, 0.88}
% Dessiner les cellules du polyomino avec leurs délimitations
% Colonne 1 (1 cellule)
\fill[lightpink] (0,0) rectangle (1,1);
\draw[black] (0,0) rectangle (1,1);
% Colonne 2 (2 cellules)
\fill[lightpink] (1,0) rectangle (2,1);
\fill[lightpink] (1,1) rectangle (2,2);
\draw[black] (1,0) rectangle (2,1);
\draw[black] (1,1) rectangle (2,2);
% Colonne 3 (3 cellules)
\fill[lightpink] (2,0) rectangle (3,1);
\fill[lightpink] (2,1) rectangle (3,2);
\fill[lightpink] (2,2) rectangle (3,3);
\draw[black] (2,0) rectangle (3,1);
\draw[black] (2,1) rectangle (3,2);
\draw[black] (2,2) rectangle (3,3);
% Colonne 4 (1 cellule)
\fill[lightpink] (3,0) rectangle (4,1);
\fill[lightpink] (3,1) rectangle (4,2);
\fill[lightpink] (3,2) rectangle (4,3);
\fill[lightpink] (3,3) rectangle (4,4);
\draw[black] (3,0) rectangle (4,1);
\draw[black] (3,1) rectangle (4,2);
\draw[black] (3,2) rectangle (4,3);
\draw[black] (3,3) rectangle (4,4);
\node at (2,-0.45){\Large 0~1~2~3};
\fill[black] (2,1) circle (5pt);
\fill[black] (3,1) circle (5pt);\fill[black] (3,2) circle (5pt);
\end{tikzpicture}}
\caption{Interior points of the $9$ polyominoes associated with $(\geq,\geq)$-Catalan words of length $4$. }
\label{fig8}
\end{figure}

In this section, we study the statistic of the number of interior points on $(\geq,\geq)$-polyominoes, i.e. those associated with Catalan words avoiding the pattern $(\geq,\geq)$. We refer to Figure~\ref{fig8} for the illustration of interior points on $(\geq,\geq)$-polyominoes of length 4. As in the previous section, for all $0\leq i \leq n-1$, $\mathcal{C}_{(\geq, \geq)}(n,i)$ is the set of the Catalan words avoiding the pattern $(\geq,\geq)$ of length $n$ whose last symbol is $i$, and we define the generating functions
\[
C_i^{\geq}(x;q)=\sum_{n\geq 1}x^n \sum_{w\in \mathcal{C}_{(\geq, \geq)}(n,i)}q^{\inter(w)},
\]
and 
\[
C^{\geq}(x;q;v)=\sum_{i \geq 0}C_i^{\geq}(x;q)v^{i}.
\]

\begin{theorem}
The generating function $C^{\geq}(x;q;v)$ is given by
\begin{align}\label{Ip}
&C^{\geq}(x;q;v)=x+x^2+\frac{x^2}{1-qv}C^{\geq}(x;q;q)+xvC^{\geq}(x;q;qv)-\frac{x^2q^2v^2}{1-qv}C^{\geq}(x;q;q^2 v).
\end{align}
\end{theorem}
\begin{proof}
According to the decomposition given in Figures \ref{f14} and \ref{f24}, we obtain the relations
\begin{align*}
&C_0^{\geq}(x;q)=x+x^2 +x^2\sum_{k \geq 1}C_{k-1}^{\geq}(x;q)q^{k-1},\\
&C_i^{\geq}(x;q)=xq^{i-1}C_{i-1}^{\geq}(x;q)+x^2q^{2i-1}C_{i-1}^{\geq}(x;q)+x^2q^{i}\sum_{k>i}C_{k-1}^{\geq}(x;q)q^{k-1}.
\end{align*}
By multiplying the last equation by $v^i$ and summing over $i\geq 1$, we obtain the functional equation
\begin{align*}
&C^{\geq}(x;q;v)-C_0^{\geq}(x;q)=xvC^{\geq}(x;q;qv)+x^2qvC^{\geq}(x;q;q^2v)\\
&\quad \quad \quad \quad \quad \quad \quad +\frac{x^2qv}{1-qv}\left(C^{\geq}(x;q;q)-C_{0}^{\geq}(x;q)\right)-\frac{x^2qv}{1-qv}\left(C^{\geq}(x;q;q^2v)-C_{0}^{\geq}(x;q)\right).
\end{align*}
Simplifying this expression we obtain the equation
\begin{align*}
&C^{\geq}(x;q;v)=x+x^2+\frac{x^2}{1-qv}C^{\geq}(x;q;q)+xvC^{\geq}(x;q;qv)-\frac{x^2q^2v^2}{1-qv}C^{\geq}(x;q;q^2 v).
\end{align*}
\end{proof}

Let $\mathcal{B}$ be the set consisting of $(\geq,\geq)$-Catalan words $w$ such that $w$ does not end with $ab$ with $a\geq b$. Then we have the following.

\begin{theorem}\label{thp} The generating function $H(x,q)$ according to the length and the number of interior points for the number of nonempty Catalan words $w$ avoiding $(\geq,\geq)$  such that $w$ does not end with $ab$, $a\geq b$ is given by 
\begin{align*}
H(x,q)=\frac{\sum_{j\geq1}x^j\prod_{i=1}^{j-1}\left(q^{i-1}-\frac{1}{1-q^i}\right)}
{1-\sum_{j\geq1}\frac{x^j}{1-q^j}\prod_{i=1}^{j-1}\left(q^{i-1}-\frac{1}{1-q^i}\right)}.
\end{align*}
\end{theorem}
\begin{proof} The proof is obtained {\it mutatis mutandis} as for Theorem~\ref{thb}, since the bijection $\psi$ (defined in the proof)  transports 
also the number of interior points.    
\end{proof}

The first terms of the series expansion of $H(x,q)$ are
\begin{align*}
x+x^2+(q+1)x^3+(q^3+q^2+q+1)x^4+(q^6+q^5+q^4+2q^3+q^2+2q+1)x^5+O(x^6).
\end{align*}

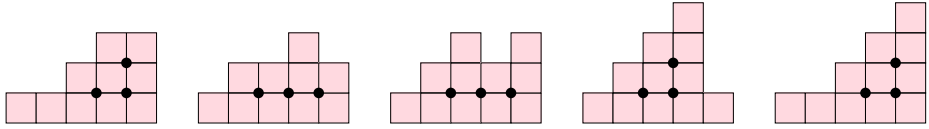
\begin{figure}[htb]
\centering
\scalebox{0.4}{\begin{tikzpicture}
% Définir la couleur de remplissage
%\definecolor{lightpink}{rgb}{1.0, 0.71, 0.76}
\definecolor{lightpink}{rgb}{1.0, 0.85, 0.88}
% Dessiner les cellules du polyomino avec leurs délimitations
% Colonne 1 (1 cellule)
\fill[lightpink] (0,0) rectangle (1,1);
\draw[black] (0,0) rectangle (1,1);
% Colonne 2 (2 cellules)
\fill[lightpink] (1,0) rectangle (2,1);
%\fill[lightpink] (1,1) rectangle (2,2);
\draw[black] (1,0) rectangle (2,1);
%\draw[black] (1,1) rectangle (2,2);
% Colonne 3 (3 cellules)
\fill[lightpink] (2,0) rectangle (3,1);
\fill[lightpink] (2,1) rectangle (3,2);
%\fill[lightpink] (2,2) rectangle (3,3);
\draw[black] (2,0) rectangle (3,1);
\draw[black] (2,1) rectangle (3,2);
%\draw[black] (2,2) rectangle (3,3);
% Colonne 4 (1 cellule)
\fill[lightpink] (3,0) rectangle (4,1);
\fill[lightpink] (3,1) rectangle (4,2);
\fill[lightpink] (3,2) rectangle (4,3);
\draw[black] (3,0) rectangle (4,1);
\draw[black] (3,1) rectangle (4,2);
\draw[black] (3,2) rectangle (4,3);
% Colonne 5 (1 cellule)
\fill[lightpink] (4,0) rectangle (5,1);
\fill[lightpink] (4,1) rectangle (5,2);
\fill[lightpink] (4,2) rectangle (5,3);
\draw[black] (4,0) rectangle (5,1);
\draw[black] (4,1) rectangle (5,2);
\draw[black] (4,2) rectangle (5,3);
\fill[black] (3,1) circle (5pt);
\fill[black] (4,1) circle (5pt);
\fill[black] (4,2) circle (5pt);
\end{tikzpicture}}
\quad 
\scalebox{0.4}{\begin{tikzpicture}
% Définir la couleur de remplissage
%\definecolor{lightpink}{rgb}{1.0, 0.71, 0.76}
\definecolor{lightpink}{rgb}{1.0, 0.85, 0.88}
% Dessiner les cellules du polyomino avec leurs délimitations
% Colonne 1 (1 cellule)
\fill[lightpink] (0,0) rectangle (1,1);
\draw[black] (0,0) rectangle (1,1);
% Colonne 2 (2 cellules)
\fill[lightpink] (1,0) rectangle (2,1);
\fill[lightpink] (1,1) rectangle (2,2);
\draw[black] (1,0) rectangle (2,1);
\draw[black] (1,1) rectangle (2,2);
% Colonne 3 (3 cellules)
\fill[lightpink] (2,0) rectangle (3,1);
\fill[lightpink] (2,1) rectangle (3,2);
%\fill[lightpink] (2,2) rectangle (3,3);
\draw[black] (2,0) rectangle (3,1);
\draw[black] (2,1) rectangle (3,2);
%\draw[black] (2,2) rectangle (3,3);
% Colonne 4 (1 cellule)
\fill[lightpink] (3,0) rectangle (4,1);
\fill[lightpink] (3,1) rectangle (4,2);
\fill[lightpink] (3,2) rectangle (4,3);
\draw[black] (3,0) rectangle (4,1);
\draw[black] (3,1) rectangle (4,2);
\draw[black] (3,2) rectangle (4,3);
% Colonne 4 (1 cellule)
\fill[lightpink] (4,0) rectangle (5,1);
\draw[black] (4,0) rectangle (5,1);
\fill[lightpink] (4,1) rectangle (5,2);
\draw[black] (4,1) rectangle (5,2);
\fill[black] (2,1) circle (5pt);
\fill[black] (3,1) circle (5pt);
\fill[black] (4,1) circle (5pt);
\end{tikzpicture}}
\quad 
\scalebox{0.4}{\begin{tikzpicture}
% Définir la couleur de remplissage
%\definecolor{lightpink}{rgb}{1.0, 0.71, 0.76}
\definecolor{lightpink}{rgb}{1.0, 0.85, 0.88}
% Dessiner les cellules du polyomino avec leurs délimitations
% Colonne 1 (1 cellule)
\fill[lightpink] (0,0) rectangle (1,1);
\draw[black] (0,0) rectangle (1,1);
% Colonne 2 (2 cellules)
\fill[lightpink] (1,0) rectangle (2,1);
\fill[lightpink] (1,1) rectangle (2,2);
\draw[black] (1,0) rectangle (2,1);
\draw[black] (1,1) rectangle (2,2);
% Colonne 3 (3 cellules)
\fill[lightpink] (2,0) rectangle (3,1);
\fill[lightpink] (2,1) rectangle (3,2);
\fill[lightpink] (2,2) rectangle (3,3);
\draw[black] (2,0) rectangle (3,1);
\draw[black] (2,1) rectangle (3,2);
\draw[black] (2,2) rectangle (3,3);
% Colonne 4 (1 cellule)
\fill[lightpink] (3,0) rectangle (4,1);
\fill[lightpink] (3,1) rectangle (4,2);
\draw[black] (3,0) rectangle (4,1);
\draw[black] (3,1) rectangle (4,2);
% Colonne 4 (1 cellule)
\fill[lightpink] (4,0) rectangle (5,1);
\draw[black] (4,0) rectangle (5,1);
\fill[lightpink] (4,1) rectangle (5,2);
\draw[black] (4,1) rectangle (5,2);
\fill[lightpink] (4,2) rectangle (5,3);
\draw[black] (4,2) rectangle (5,3);
\fill[black] (2,1) circle (5pt);
\fill[black] (3,1) circle (5pt);
\fill[black] (4,1) circle (5pt);
\end{tikzpicture}}
\quad\scalebox{0.4}{\begin{tikzpicture}
% Définir la couleur de remplissage
%\definecolor{lightpink}{rgb}{1.0, 0.71, 0.76}
\definecolor{lightpink}{rgb}{1.0, 0.85, 0.88}
% Dessiner les cellules du polyomino avec leurs délimitations
% Colonne 1 (1 cellule)
\fill[lightpink] (0,0) rectangle (1,1);
\draw[black] (0,0) rectangle (1,1);
% Colonne 2 (2 cellules)
\fill[lightpink] (1,0) rectangle (2,1);
\fill[lightpink] (1,1) rectangle (2,2);
\draw[black] (1,0) rectangle (2,1);
\draw[black] (1,1) rectangle (2,2);
% Colonne 3 (3 cellules)
\fill[lightpink] (2,0) rectangle (3,1);
\fill[lightpink] (2,1) rectangle (3,2);
\fill[lightpink] (2,2) rectangle (3,3);
\draw[black] (2,0) rectangle (3,1);
\draw[black] (2,1) rectangle (3,2);
\draw[black] (2,2) rectangle (3,3);
% Colonne 4 (1 cellule)
\fill[lightpink] (3,0) rectangle (4,1);
\fill[lightpink] (3,1) rectangle (4,2);
\fill[lightpink] (3,2) rectangle (4,3);
\fill[lightpink] (3,3) rectangle (4,4);
\draw[black] (3,0) rectangle (4,1);
\draw[black] (3,1) rectangle (4,2);
\draw[black] (3,2) rectangle (4,3);
\draw[black] (3,3) rectangle (4,4);
% Colonne 4 (1 cellule)
\fill[lightpink] (4,0) rectangle (5,1);
\draw[black] (4,0) rectangle (5,1);
\fill[black] (2,1) circle (5pt);
\fill[black] (3,1) circle (5pt);
\fill[black] (3,2) circle (5pt);
\end{tikzpicture}}
\quad \scalebox{0.4}{\begin{tikzpicture}
% Définir la couleur de remplissage
%\definecolor{lightpink}{rgb}{1.0, 0.71, 0.76}
\definecolor{lightpink}{rgb}{1.0, 0.85, 0.88}
% Dessiner les cellules du polyomino avec leurs délimitations
% Colonne 1 (1 cellule)
\fill[lightpink] (0,0) rectangle (1,1);
\draw[black] (0,0) rectangle (1,1);
% Colonne 2 (2 cellules)
\fill[lightpink] (1,0) rectangle (2,1);

\draw[black] (1,0) rectangle (2,1);

% Colonne 3 (3 cellules)
\fill[lightpink] (2,0) rectangle (3,1);
\draw[black] (2,0) rectangle (3,1);
\fill[lightpink] (2,1) rectangle (3,2);
\draw[black] (2,1) rectangle (3,2);

% Colonne 4 (1 cellule)
\fill[lightpink] (3,0) rectangle (4,1);
\fill[lightpink] (3,1) rectangle (4,2);
\fill[lightpink] (3,2) rectangle (4,3);
\draw[black] (3,0) rectangle (4,1);
\draw[black] (3,1) rectangle (4,2);
\draw[black] (3,2) rectangle (4,3);
% Colonne 5 (1 cellule)
\fill[lightpink] (4,0) rectangle (5,1);
\fill[lightpink] (4,1) rectangle (5,2);
\fill[lightpink] (4,2) rectangle (5,3);
\fill[lightpink] (4,3) rectangle (5,4);
\draw[black] (4,0) rectangle (5,1);
\draw[black] (4,1) rectangle (5,2);
\draw[black] (4,2) rectangle (5,3);
\draw[black] (4,3) rectangle (5,4);
\fill[black] (3,1) circle (5pt);
\fill[black] (4,1) circle (5pt);
\fill[black] (4,2) circle (5pt);
\end{tikzpicture}}
\caption{The $5$ polyominoes of length $5$ with 3 interior points . }\label{fig9}
\end{figure}

\begin{theorem}\label{thmintp}The generating function for the number of nonempty $(\geq,\geq)$-Catalan words according to the length and the number of interior points is given by $$C^{\geq}(x,q;1)=  \sum_{i=1}^{\infty} x^i q^{\frac{(i-2)(i-1)}{2}}  \prod_{j=0}^{i-1} \left(1+\mathit{H} \left(x q^j, q \right)\right).
$$
 \end{theorem}
  \begin{proof} Any nonempty Catalan word avoiding $(\geq,\geq)$ is of the form either $u0$ or $u0(1+v)$ where $u\in\mathcal{B}$ and $v\in \mathcal{C}$, which generate the functional equation  
$$C^{\geq}(x,q;1)=x(1+H(x,q))+\frac{x}{q}C^{\geq}(xq,q;1)(1+H(x,q)).$$
Iterating this equation we obtain the desired result.
 \end{proof}
The first terms of the series expansion of $C^{\geq}(x,q;1)$ are
\begin{align*}
x+2x^2+(2q+2)x^3+(2q^3+2q^2+3q+2)x^4+(2q^6+2q^5+3q^4+\mathbf{5q^3}+3q^2+4q+2)x^5+\\(2q^{10}+2q^9+3q^8+5q^7+6q^6+7q^5+6q^4+7q^3+6q^2+5q+2)x^6+O(x^7).
\end{align*}
We refer to Figure~\ref{fig9} for an illustration of the $(\geq,\geq)$-polyominoes of length 5 with 3 interior points.

Let $p(n, i)$ denote the total number of interior points of the $(\geq,\geq)$-polyominoes of length $n$ that end with a column of height $i$. The first
few values are
\[
[p(n,i)]_{n,i\geq 1}=\left[
                            \begin{array}{cccccc}
                              0 & 0 & 0 & 0 & 0  &\cdots \\
                              0 & 0 & 0 & 0 & 0 &\cdots\\
                              0 & 1 & 1 & 0 & 0 &\cdots\\
                              1 & 3 & 6 & 3 & 0 &\cdots \\
                              6& 11 & 18 & 18 & 6&\cdots\\
                              \vdots & \vdots & \vdots & \vdots & \vdots & \ddots
                            \end{array}
                          \right].
\]

From the decomposition given in Figures~\ref{f14} and \ref{f24}, we have for $n\geq2$ and $2 \leq i \leq n$,
\begin{align*}
&p(n,i)=p(n-1,i-1)+(i-2)c(n-1,i-2)+\sum_{k=i}^{n-2}\left(p(n-2,k-1)+(i+k-3)c(n-2,k-2)\right).
\end{align*}

If we consider the difference $p(n,i)-p(n,i-1)$, then for $n\geq2$ and $3\leq i\leq n$, we obtain
the recurrence relation

\begin{align*}
&p(n,i)=p(n,i-1)+p(n-1,i-1)-p(n-1,i-2)+(i-2)c(n-1,i-2)\\
&\quad \quad-(i-3)c(n-1,i-3)+p(n-2,n-3)-p(n-2,i-2)+(n+i-5)c(n-2,n-4)\\
&\quad \quad -(2i-4)c(n-2,i-3)+\sum_{k=i-2}^{n-2}c(n-2,k-2).
\end{align*}

Let $p(n)$ be the total number of interior points of over all $(\geq,\geq)$-polyominoes of length~$n$. The first few values for $1\leq n\leq 10$ are
\[
0,\quad 0,\quad 2,\quad 13,\quad 59,\quad 230,\quad 830, \quad 2858, \quad 9547, \quad 31227.
\]
The following theorem gives a combinatorial formula to calculate the sequence $p(n)$.
\begin{theorem} The generating function for the total number of interior points over all $(\geq,\geq)$-Catalan words with respect to the length is 
\begin{align*}
P(x,1)=\frac{\left(-8 x^{4}-8 x^{3}+5 x^{2}+3 x -1\right) \sqrt{-3 x^{2}-2 x +1}-6 x^{6}-10 x^{5}+12 x^{4}+17 x^{3}-4 x^{2}-4 x +1}{6 x^{6}+4 x^{5}-2 x^{4}},\end{align*}
and the coefficient of $x^n$ is given by
\[
 p(n)=\frac{1}{2}\left(3^{n+1}+8T_n+8T_{n+1}-5T_{n+2}-3T_{n+3}+T_{n+4}\right).
\]
An asymptotic approximation is $3^{n+1}/2.$
\end{theorem}

\begin{proof}
Let $P(x,v)=\left.\frac{\partial}{\partial q}C^{\geq}(x;q;v) \right|_{q=1}$. Then by differentiating \eqref{Ip} with respect to $q$, we obtain
\begin{align*}
P(x;v)&=\frac{x^2v}{(1-v)^2}(M(x)-1)+\frac{x^2}{1-v}\left(P(x,1)+\left.\frac{\partial}{\partial v}C^{\geq}(x;1;v) \right|_{v=1}\right)\\
&\quad +xv\left(P(x,v)+v\frac{\partial}{\partial v}C^{\geq}(x;1;v) \right)-\frac{x^2v^2}{1-v}\left(P(x,v)+2v\frac{\partial}{\partial v}C^{\geq}(x;1;v) \right)\\
&\quad +\frac{x^2v^2(v-2)}{(1-v)^2}C^{\geq}(x;1;v).
\end{align*}
Notice that with relation \eqref{gf3}, we obtain
\begin{align*}
&\frac{\partial}{\partial v} C^{\geq}(x;1;v)=\frac{\partial}{\partial v}\left(\frac{x(1-v)-x^2v+x^2M(x)}{1-v-xv(1-v)+x^2v^2}\right) \label{f1}\\ 
&\quad =\frac{\left(x +1\right) \left(\left(x v -\frac{1}{2}\right) \sqrt{-3 x^{2}-2 x +1}+\frac{1}{2}+v^{2} x^{3}+\left(v^{2}-v \right) x^{2}+\left(-v -\frac{1}{2}\right) x \right)}{\left(v^{2} x^{2}+\left(v^{2}-v \right) x -v +1\right)^{2}}.
\end{align*}

Now, we group the terms $P(x;v)$ and we rewrite the equation so that the left-hand side is 
$$\frac{P \! \left(x , v\right) \left(1-v -x v \left(-v +1\right)+v^{2} x^{2}\right)^3}{-v +1}.$$

By twice differentiating this equation with respect to $v$ and canceling the factor $1-v-xv(1-v)+x^2v^2$ with $$v={\frac {1+x-\sqrt {-3\,{x}^{2}-2\,x+1}}{2x \left( x+1 \right) }}=\frac{1+xM(x)}{1+x},$$ we obtain that
\begin{align*}
P(x,1)&=\frac{\left(-8 x^{4}-8 x^{3}+5 x^{2}+3 x -1\right) \sqrt{-3 x^{2}-2 x +1}-6 x^{6}-10 x^{5}+12 x^{4}+17 x^{3}-4 x^{2}-4 x +1}{6 x^{6}+4 x^{5}-2 x^{4}}
. \end{align*}
The following decomposition
$$P(x,1)=\frac{-6 x^{5}-4 x^{4}+16 x^{3}+x^{2}-5 x +1}{6 x^{5}-2 x^{4}}+\frac{8 x^{4}+8 x^{3}-5 x^{2}-3 x +1}{2 x^{4}}T(x)$$
allows to obtain the desired result (see Section \ref{sec:2} for the definition of $T(x)$).
\end{proof} 
 
\section*{Acknowledgment}
The first and third authors would like to thank the Laboratoire d’Informatique de Bourgogne for the warm hospitality during their visit at Université de Bourgogne where part of this work was done. This work was supported in part by DG-RSDT (Algeria), PRFU Project, No. C00L03UN180120220002 and by the project ANR PiCs ANR-22-CE48-0002 (France).


\begin{thebibliography}{20}

\bibitem{ArcBleKnop} M.~Archibald, A.~Blecher, and A.~Knopfmacher. Parameters in inversion sequences, Math. Slovaca \textbf{73} (3) (2023), 551--564.

\bibitem{Kernel} C.~Banderier, M.~Bousquet-M\'elou, A.~Denise, P.~Flajolet, D.~Gardy, and D.~Gouyou-Beauchamps. Generating functions for generating trees, Discrete Math. \textbf{246} (2002), 29--55.

\bibitem{Banfla}  C.~Banderier and P.~Flajolet. Basic analytic combinatorics of directed lattice paths, Theoret. Comput. Sci. \textbf{281} (2002), 37--80.

\bibitem{BKR}  J.-L.~Baril, S. Kirgizov, J.~L.~Ram\'irez and D. Villamizar. The combinatorics of Motzkin polyominoes, Discrete Applied Math. \textbf{364} (2025), 1--15.


\bibitem{relation} J.-L.~Baril and J.~L.~Ram\'irez. Descent distribution on Catalan words avoiding ordered pairs of relations, Adv. in Applied Math. \textbf{149} (2023), 102551.
 
\bibitem{RiordanN} F.~R.~Bernhart. Catalan, Motzkin, and Riordan numbers, Discrete Math. \textbf{204} (1999), 73--112.

\bibitem{BLE3} A.~Blecher, C.~Brennan, and A.~Knopfmacher. Combinatorial parameters in bargraphs, Quaest. Math. \textbf{39} (2016), 619--635.

\bibitem{BleBreKnop}  A.~Blecher, C.~Brennan, and A.~Knopfmacher. The site-perimeter of compositions, Discrete Math. Appl. \textbf{32} (2) (2022), 75--89.

\bibitem{BleBreKnop3} A.~Blecher, C.~Brennan, A.~Knopfmacher, and T.~Mansour.  The perimeter of words, Discrete Math. \textbf{340} (10) (2017), 2456--2465. 

\bibitem{BlKn} A.~Blecher and A.~~Knopfmacher. Cells of fixed height in Catalan words and restricted growth functions, Advances in Applied Mathematics, \textbf{164} (2025), 102835.

\bibitem{Bou} M.~Bousquet-Mélou and A.~Rechnitzer. The site-perimeter of bargraphs,
Adv. in Appl. Math. \textbf{31}(2003), 86--112.

\bibitem{CallManRam} D.~Callan, T.~Mansour, and J.~L.~Ram\'irez.  Statistics on bargraphs of Catalan words, J. Autom. Lang. Comb. \textbf{26} (2021), 177--196.

\bibitem{Don} R.~Donaghey and L.~W.~Shapiro. Motzkin numbers, J. Combin. Theory Ser. A  \textbf{23}(3) (1977),  291--301.

\bibitem{Fla} P.~Flajolet and R.~Sedgewick. \emph{Analytic Combinatorics}, Cambridge University Press, 2009.

\bibitem{Ges} I.~M.~Gessel and S.~Ree. Lattice paths and Faber polynomials,  Advances in Combinatorial Methods and Applications to Probability and Statistics, Birkhauser Verlag, Boston, 1997.


 \bibitem{Book1} A.~J.~Guttmann (Ed.) \emph{Polygons, Polyominoes and Polycubes}, Lecture Notes in Physics 775. Springer, Heidelberg, Germany, 2009.

\bibitem{Krat} C.~Krattenthaler.
Lattice path enumeration, 
 In: B\'ona, M. (ed.) Handbook of Enumerative Combinatorics, CRC Press, (2015).


\bibitem{ManA} T.~Mansour.   Semi-perimeter and inner site-perimeter of $k$-ary words and bargraphs, Art Disc. Appl. Math. \textbf{4} (2021), Article P1.06.

\bibitem{ManA2} T.~Mansour.  The perimeter and the site-perimeter of set partitions, Electron. J. Comb. \textbf{26} (2) (2019), Article \#P2.30.

\bibitem{ManA3} T.~Mansour.  Interior vertices in set partitions, Adv. in Applied Math. \textbf{101} (2018), 60--69.

  
\bibitem{ManSha2} T.~Mansour and A.~Sh.~Shabani. Enumerations on bargraphs, Discrete Math. Lett. \textbf{2} (2019), 65--94.

\bibitem{ManRam} T.~Mansour and J.~L.~Ram\'irez. Enumerations on polyominoes determined by Fuss-Catalan words, Australas. J. Comb. \textbf{81} (2021), 447--457.

 
\bibitem{ManRamMot} T.~Mansour and J.~L. Ram\'irez. Exterior corners on bargraphs  of Motzkin words,  In: F.~Hoffman,  S.~Holliday, Z.~Rosen, F.~Shahrokhi, and J.~Wierman (eds). Combinatorics, Graph Theory and Computing 2021. Springer Proceedings in Mathematics \& Statistics, vol 448. Springer, (2024).

 \bibitem{Toc}
 T.~Mansour,  J.~L.~Ram\'irez, and D.~A. Toquica. Counting lattice   points on bargraphs of {C}atalan words,  Math. Comput. Sci. \textbf{15} (2021), 701--713.


  
\bibitem{Orl} A.~G.~Orlov. On asymptotic behavior of the Taylor coefficients of algebraic functions,  Sib. Math. J.
\textbf{25}(5) (1994), 1002--1013.



\bibitem{pro} H.~Prodinger.
\newblock The kernel method: a collection of examples, 
 Sém. Lothar. Combin. \textbf{50} (2004), Paper B50f.


\bibitem{OEIS} N.~J.~A.~Sloane. The On-Line Encyclopedia of Integer Sequences, \url{http://oeis.org/}.

\bibitem{stan} R.~P.~Stanley. \emph{Catalan Numbers}, Cambridge University Press, 2015.




\end{thebibliography}
\end{document}